\documentclass[secthm]{elsart}
\usepackage{amssymb}
\usepackage{amsfonts}
\usepackage{amsmath}

\input{tcilatex}
\numberwithin{equation}{section}

\begin{document}

\begin{frontmatter}
\title{Boolean Inner-Product Spaces and Boolean Matrices}
\author{Stan Gudder}
\address{Department of Mathematics, University of Denver, Denver CO 80208}
\ead{sgudder@math.du.edu}
\author{Fr\'ed\'eric Latr\'emoli\`ere}
\address{Department of Mathematics, University of Denver, Denver CO 80208}
\ead{frederic@math.du.edu}
\begin{abstract}
This article discusses the concept of Boolean spaces endowed with a Boolean
valued inner product and their matrices. A natural inner product structure
for the space of Boolean $n$-tuples is introduced. Stochastic boolean
vectors and stochastic and unitary Boolean matrices are studied. A dimension
theorem for orthonormal bases of a Boolean space is proven. We characterize
the invariant stochastic Boolean vectors for a Boolean stochastic matrix and
show that they can be used to reduce a unitary matrix. Finally, we obtain a
result on powers of stochastic and unitary matrices.
\end{abstract}
\end{frontmatter}

\section{Introduction}

\qquad A Boolean space $\mathcal{L}_{n}\left( \mathcal{B}\right) $ is the
set of all $n$-tuples of elements of a fixed Boolean algebra $\mathcal{B}$.
The elements of $\mathcal{L}_{n}\left( \mathcal{B}\right) $ are called
Boolean vectors and they possess a natural linear space-like structure.
Moreover, we can define on $\mathcal{L}_{n}\left( \mathcal{B}\right) $ an
operation which is analogous to an inner product. By using this
\textquotedblleft inner product\textquotedblright\ we can also define a $%
\mathcal{B}$-valued norm and orthogonality relations for Boolean vectors.

\qquad A Boolean matrix is a matrix whose entries are elements of a Boolean
algebra $\mathcal{B}$. With the natural choice of matrix multiplication
defined in terms of the lattice operations of $\mathcal{B}$, such matrices
become the linear mappings between Boolean linear spaces. The study of
Boolean matrices is a fascinating blend of linear algebra and boolean
algebra which finds many applications, and was undertaken in \cite%
{Blyth67,Cechlarova03,Giveon64,Jagannadham66,Luce52,Rutherford63,Rutherford63b,Rutherford64,Sindak75,Skornyakov86,Subrahmanyam64,Subrahmanyam65,Subrahmanyam67,Wedderburn34,Tan98,Tan01,Yoeli61,Gregory94}
.

\qquad An important concept in our work is that of a stochastic vector.
These are Boolean vectors of norm one whose components are mutually
disjoint. In particular, a finite partition of the universe of a Boolean
algebra would correspond to a stochastic Boolean vector. We define an
orthonormal basis of $\mathcal{L}_{n}\left( \mathcal{B}\right) $ the usual
way and it turns out it must be made of stochastic vectors. Our first main
result is that all orthonormal bases for $\mathcal{L}_{n}\left( \mathcal{B}
\right) $ have cardinality $n$ and conversely, any orthonormal set of
stochastic vectors with cardinality $n$ is a basis for $\mathcal{L}
_{n}\left( \mathcal{B}\right) $. Our next main result states that any
orthonormal set of stochastic vectors in $\mathcal{L}_{n}\left( \mathcal{B}
\right) $ can be extended to an orthonormal basis for $\mathcal{L}_{n}\left( 
\mathcal{B}\right) $. In order to prove this result, we introduce a notion
of linear subspace of $\mathcal{L}_{n}\left( \mathcal{B}\right) $.

\qquad We define stochastic and unitary Boolean matrices in terms of
properties of their product with their adjoint matrices. We then show that
stochastic Boolean matrices are precisely those whose columns are stochastic
vectors and unitary matrices are precisely those whose rows and columns are
stochastic.

\qquad We next characterize the invariant stochastic Boolean vectors for
stochastic Boolean matrices and show that they can be employed to reduce
unitary Boolean matrices. As mentioned in Section 2, stochastic Boolean
matrices may be used to describe a dynamics analogous to a Markov chain. It
is thus of interest to consider powers of stochastic Boolean matrices
because they correspond to iterations in the dynamics. Our last result
concerns such powers. The paper includes examples that illustrate various
points which we wish to emphasize.

\qquad As a matter of notations, we shall write $\mathbb{N}$ as the set of
nonzero natural numbers.

\section{Definitions and Motivation}

\qquad Throughout this article, $\mathcal{B}$ will denote a Boolean algebra.
We denote the smallest and largest element of $\mathcal{B}$ respectively by $%
0$ and $1$. For any $a\in \mathcal{B}$, we denote by $a^{c}$ its complement.
For $a,b\in \mathcal{B}$, we denote the infimum of $a$ and $b$ by $ab$
(instead of $a\wedge b$). We denote by $a\backslash b=a\left( b^{c}\right) $
. The supremum of $a,b$ is denoted by $a\vee b$.

\qquad For all $n\in \mathbb{N}$ we denote by $\mathcal{L}_{n}\left( 
\mathcal{B}\right) $ the set of all $n$-tuples of elements in $\mathcal{B}$.
We endow $\mathcal{L}_{n}\left( \mathcal{B}\right) $ with the following
operations: if $\underline{a}=\left( a_{1},\ldots ,a_{n}\right) $ and $%
\underline{b}=\left( b_{1},\ldots ,b_{n}\right) $ are in $\mathcal{L}
_{n}\left( \mathcal{B}\right) ,$ and $c\in \mathcal{B}$ then 
\begin{equation*}
\underline{a}+\underline{b}=\left( a_{1}\vee b_{1},\ldots ,a_{n}\vee
b_{n}\right)
\end{equation*}
and 
\begin{equation*}
c\underline{a}=\left( ca_{1},\ldots ,ca_{n}\right) \text{.}
\end{equation*}
Then $\mathcal{L}_{n}\left( \mathcal{B}\right) $ has the usual properties of
a linear space except for the lack of additive inverses. In particular, our
structure differs from the notion of Boolean vector space introduced in \cite%
{Subrahmanyam64,Subrahmanyam65,Subrahmanyam67} which assumes an underlying
additive group and is best modelled by the action of a Boolean space on a
regular vector space by means of a (finitely additive) measure.

\qquad We call the elements of $\mathcal{L}_{n}\left( \mathcal{B}\right) $ 
\emph{Boolean vectors} and call $\mathcal{L}_{n}\left( \mathcal{B}\right) $
a \emph{Boolean (linear) space}. We will use the following definitions
throughout this paper

\begin{defn}
A Boolean vector $\underline{a}=\left( a_{1},\ldots ,a_{n}\right) $ is an
orthovector when $a_{i}a_{j}=0$ for $i,j\in \left\{ 1,\ldots ,n\right\} $
and $i\not=j$.
\end{defn}

\begin{defn}
An orthovector $\underline{a}=\left( a_{1},\ldots ,a_{n}\right) $ is a
stochastic vector when $\dbigvee\limits_{i=1}^{n}a_{i}=1$.
\end{defn}

\qquad The Boolean space $\mathcal{L}_{n}\left( \mathcal{B}\right) $ is
endowed with a natural inner product.

\begin{defn}
Let $\underline{a}=\left( a_{1},\ldots ,a_{n}\right) $ and $\underline{b}
=\left( b_{1},\ldots ,b_{n}\right) $ in $\mathcal{L}_{n}\left( \mathcal{B}
\right) $. Then we define the $\mathcal{B}$-valued inner product of these
two vectors by 
\begin{equation*}
\left\langle \underline{a},\underline{b}\right\rangle
=\dbigvee\limits_{i=1}^{n}a_{i}b_{i}\text{.}
\end{equation*}

The norm of \underline{$a$} is defined by $\left\Vert \underline{a}
\right\Vert =\left\langle \underline{a},\underline{a}\right\rangle $.
\end{defn}

\bigskip \qquad The Boolean inner product shares most of the usual
properties of the Euclidian inner product, if we replace scalar sums and
products by the supremum and infimum in $\mathcal{B}$. Thus given $%
\underline{a},\underline{b},\underline{c}\in \mathcal{L}_{n}\left( \mathcal{%
B }\right) $ and $\alpha \in \mathcal{B}$ then

\begin{itemize}
\item $\left\langle \alpha \underline{a}+\underline{b},\underline{c}
\right\rangle =\alpha \left\langle \underline{a},\underline{c}\right\rangle
\vee \left\langle \underline{b},\underline{c}\right\rangle $,

\item $\left\langle \underline{a},\underline{b}\right\rangle =\left\langle 
\underline{b},\underline{a}\right\rangle $,

\item $\left\langle \alpha \underline{a},\underline{c}\right\rangle
=\left\langle \underline{a},\alpha \underline{c}\right\rangle $,

\item $\left\langle \underline{a},\underline{a}\right\rangle =0$ if and only
if $\underline{a}=\left( 0,\ldots ,0\right) =\underline{0}$.
\end{itemize}

We now give some properties of the norm.

\begin{thm}
\label{Norm}Let $\underline{a},$ $\underline{b}\in \mathcal{L}_{n}\left( 
\mathcal{B}\right) $ and $c\in \mathcal{B}$. Then

\begin{enumerate}
\item $\left\Vert c\underline{a}\right\Vert =c\left\Vert \underline{a}
\right\Vert $,

\item $\left\Vert \underline{a}+\underline{b}\right\Vert =\left\Vert 
\underline{a}\right\Vert \vee \left\Vert \underline{b}\right\Vert $,

\item $\left\langle \underline{a},\underline{b}\right\rangle \leq \left\Vert 
\underline{a}\right\Vert \left\Vert \underline{b}\right\Vert $,

\item If $\underline{a}$ and $\underline{b}$ are orthovectors and $%
\left\Vert \underline{a}\right\Vert =\left\Vert \underline{b}\right\Vert $
then $\left\langle \underline{a},\underline{b}\right\rangle =\left\Vert 
\underline{a}\right\Vert \left\Vert \underline{b}\right\Vert $ if and only
if $\underline{a}=\underline{b}$.
\end{enumerate}
\end{thm}

\begin{pf}
We have 
\begin{equation*}
\left\Vert c\underline{a}\right\Vert =\left\langle c\underline{a},c 
\underline{a}\right\rangle =c\left\langle \underline{a},c\underline{a}
\right\rangle =c\left\langle \underline{a},\underline{a}\right\rangle
=c\left\Vert \underline{a}\right\Vert
\end{equation*}
and, denoting $\underline{a}=\left( a_{1},\ldots ,a_{n}\right) $ and $%
\underline{b}=\left( b_{1},\ldots ,b_{n}\right) $, we have 
\begin{equation*}
\left\Vert \underline{a}+\underline{b}\right\Vert
=\dbigvee\limits_{i=1}^{n}\left( a_{i}\vee b_{i}\right) =\left(
\dbigvee\limits_{i=1}^{n}a_{i}\right) \vee \left(
\dbigvee\limits_{i=1}^{n}b_{i}\right) =\left\Vert \underline{a}\right\Vert
\vee \left\Vert \underline{b}\right\Vert
\end{equation*}
while 
\begin{equation*}
\left\langle \underline{a},\underline{b}\right\rangle
=\dbigvee\limits_{i=1}^{n}a_{i}b_{i}\leq
\dbigvee\limits_{i,j=1}^{n}a_{i}b_{j}=\left(
\dbigvee\limits_{i=1}^{n}a_{i}\right) \left(
\dbigvee\limits_{j=1}^{n}b_{j}\right) =\left\Vert \underline{a}\right\Vert
\left\Vert \underline{b}\right\Vert \text{.}
\end{equation*}

Now let us assume that $\underline{a}$ and $\underline{b}$ are orthovectors
and $\left\Vert \underline{a}\right\Vert =\left\Vert \underline{b}
\right\Vert $ and that $\left\langle \underline{a},\underline{b}
\right\rangle =\left\Vert \underline{a}\right\Vert \left\Vert \underline{b}
\right\Vert $. Hence, $\left\langle \underline{a},\underline{b}\right\rangle
=\left\Vert \underline{a}\right\Vert $ so $\dbigvee
\limits_{i=1}^{n}a_{i}b_{i}=\dbigvee\limits_{i=1}^{n}a_{i}$. Hence, for all $%
j\in \left\{ 1,\ldots ,n\right\} $ 
\begin{equation*}
a_{j}b_{j}=\left( \dbigvee\limits_{i=1}^{n}a_{i}b_{i}\right)
b_{j}=a_{j}b_{j}\vee \left( \dbigvee\limits_{i\not=j}a_{i}b_{j}\right) \text{
.}
\end{equation*}
Hence $\dbigvee\limits_{i\not=j}a_{i}b_{j}\leq a_{j}b_{j}$ yet $%
\dbigvee\limits_{i\not=j}a_{i}b_{j}\leq a_{j}^{c}b_{j}$ since $\underline{a}$
is an orthovector, so $\dbigvee\limits_{i\not=j}a_{i}b_{j}=0$ and thus $%
a_{i}b_{j}=0$. Therefore 
\begin{equation*}
a_{j}\left( \left\Vert \underline{a}\right\Vert \backslash b_{j}\right)
=a_{j}\left( \left\Vert \underline{b}\right\Vert \backslash b_{j}\right)
=a_{j}\left( \dbigvee\limits_{i\not=j}b_{i}\right)
=\dbigvee\limits_{i\not=j}a_{j}b_{i}=0\text{.}
\end{equation*}
Hence, using again that $a$ is an orthovector, $a_{j}=a_{j}\left\Vert 
\underline{a}\right\Vert =a_{j}b_{j}\leq b_{j}$. Symmetrically, $b_{j}\leq
a_{j}$ so $a_{j}=b_{j}$ for all $j\in \left\{ 1,\ldots ,n\right\} $. Hence $%
\underline{a}=\underline{b}$.\qed
\end{pf}

\bigskip \qquad Note that the condition $\left\Vert \underline{a}\right\Vert
=\left\Vert \underline{b}\right\Vert $ in the last statement of Theorem (\ref%
{Norm}) is necessary. If we let $\underline{a}=\left( a,0,\ldots ,0\right) $
and $\underline{b}=\left( b,0,\ldots ,0\right) $ with $a,b\in \mathcal{B}$
and $a\not=b$ then $\underline{a},\underline{b}$ are orthovectors of
different norms, and yet trivially $\left\langle \underline{a},\underline{b}
\right\rangle =\left\Vert \underline{a}\right\Vert \left\Vert \underline{b}
\right\Vert $. Also, the condition that $\underline{a}$ and $\underline{b}$
are orthovectors is necessary since if $\underline{a}=\left( 1,a\right) $
and $\underline{b}=\left( 1,b\right) $ for $a,b\in \mathcal{B}$ with $%
a\not=b $ then $\left\Vert \underline{a}\right\Vert =\left\Vert \underline{b}
\right\Vert =1$ and $\left\langle \underline{a},\underline{b}\right\rangle
=1 $.

\begin{cor}
If $\underline{a}$ and $\underline{b}$ are stochastic Boolean vectors then $%
\left\langle \underline{a},\underline{b}\right\rangle =1$ if and only if $%
\underline{a}=\underline{b}$.
\end{cor}

\begin{pf}
By assumption, $\underline{a}$ and $\underline{b}$ are orthovectors with $%
\left\Vert \underline{a}\right\Vert =\left\Vert \underline{b}\right\Vert =1$
so the result follows from Theorem (\ref{Norm}).\qed
\end{pf}

We now introduce the following standard notions:

\begin{defn}
Two vectors $\underline{a}$ and $\underline{b}$ in $\mathcal{L}_{n}\left( 
\mathcal{B}\right) $ are orthogonal when $\left\langle \underline{a}, 
\underline{b}\right\rangle =0$, in which case we shall write $\underline{a}
\perp \underline{b}$. The vector $\underline{a}$ is a unit vector when $%
\left\Vert \underline{a}\right\Vert =1$.
\end{defn}

\begin{defn}
An orthogonal set in $\mathcal{L}_{n}\left( \mathcal{B}\right) $ is a subset 
$E$ of $\mathcal{L}_{n}\left( \mathcal{B}\right) $ such that for all $%
\underline{e},\underline{f}\in E$ we have $\underline{e}\not=\underline{f}
\implies \left\langle \underline{e},\underline{f}\right\rangle =0$. An
orthonormal subset of $\mathcal{L}_{n}\left( \mathcal{B}\right) $ is an
orthogonal set whose elements all have norm $1$.
\end{defn}

\qquad The next section of this paper will address the concept of dimension
for a Boolean vector space. It will be based on the notion of basis. We now
introduce:

\begin{defn}
Let $\mathcal{A}$ be a subset of $\mathcal{L}_{n}\left( \mathcal{B}\right) $
. A vector $\underline{b}\in \mathcal{L}_{n}\left( \mathcal{B}\right) $ is a
linear combination of elements in $\mathcal{A}$ when there exists a finite
subset $\left\{ \underline{a_{1}},\ldots ,\underline{a_{m}}\right\} $ of $%
\mathcal{A}$ and $b_{1},\ldots ,b_{m}\in \mathcal{B}$ such that $\underline{%
b }=\sum_{i=1}^{m}b_{i}\underline{a_{i}}$.

A subset $\mathcal{A}$ of $\mathcal{L}_{n}\left( \mathcal{B}\right) $ is a
generating subset of $\mathcal{L}_{n}\left( \mathcal{B}\right) $ when all
vectors in $\mathcal{L}_{n}\left( \mathcal{B}\right) $ are linear
combinations of elements in $\mathcal{A}$.

A subset $\mathcal{A}$ is free when for any $b_{i},d_{j}\in \mathcal{B}
\backslash \left\{ 0\right\} $ and $\underline{a_{i}},\underline{c_{j}}\in 
\mathcal{A}$ with $i=1,\ldots ,m$ and $j=1,\ldots k$ such that $%
\sum_{i=1}^{m}b_{i}\underline{a_{i}}=\sum_{j=1}^{k}d_{j}\underline{c_{j}}$
we have: 
\begin{equation*}
m=k\text{, }\left\{ b_{1},\ldots ,b_{m}\right\} =\left\{ d_{1},\ldots
,d_{m}\right\} \text{ and }\left\{ \underline{a_{1}},\ldots ,\underline{%
a_{m} }\right\} =\left\{ \underline{c_{1}},\ldots ,\underline{c_{m}}\right\} 
\text{ .}
\end{equation*}
\end{defn}

\qquad Thus a set $\mathcal{A}$ is free whenever a linear combination of
elements in $\mathcal{A}$ has unique nonzero coefficients and associated
vectors of $\mathcal{A}$. We naturally introduce:

\begin{defn}
\label{Basis}A subset $\mathcal{A}$ of $\mathcal{L}_{n}\left( \mathcal{B}
\right) $ is a basis of $\mathcal{L}_{n}\left( \mathcal{B}\right) $ when
every element of $\mathcal{L}_{n}\left( \mathcal{B}\right) $ can be written
as a unique linear combination of elements of $\mathcal{A}$ with nonzero
coefficients, i.e. when $\mathcal{A}$ is generating and free.
\end{defn}

\qquad A first easy observation is that a basis must be made of unit vectors.

\begin{lem}
Let $\mathcal{A}$ be a basis of $\mathcal{L}_{n}\left( \mathcal{B}\right) $.
If $\underline{a}\in \mathcal{A}$ then $\left\Vert \underline{a}\right\Vert
=1$.
\end{lem}

\begin{pf}
Note first that, if $\underline{0}=(0,\ldots ,0)\in \mathcal{L}_{n}\left( 
\mathcal{B}\right) $ were in $\mathcal{A}$ and $\underline{1}=(1,\ldots
,1)\in \mathcal{L}_{n}\left( \mathcal{B}\right) $ then $\underline{1}= 
\underline{1}=\underline{1}+\underline{0}$, so $\underline{1}$ could be
written as two distinct linear combinations of elements in $\mathcal{A}$
with coefficients $1$. This is a contradiction so $\underline{0}\not\in 
\mathcal{A}$. Let $\underline{a}\in \mathcal{A}$. Then $\underline{a}=1 
\underline{a}=\left\Vert \underline{a}\right\Vert \underline{a}$. Hence if $%
\left\Vert \underline{a}\right\Vert \not=1$ then $\underline{a}$ can be
written as two distinct linear combinations of elements in $\mathcal{A}$
with nonzero coefficients (since $\underline{a}\not=\underline{0}$ so $%
\left\Vert \underline{a}\right\Vert \not=0$) which contradicts the
definition of a basis.\qed
\end{pf}

\bigskip A second easy observation is:

\begin{lem}
\label{OrthoFree}Let $\mathcal{A}$ be an orthonormal set in $\mathcal{L}
_{n}\left( \mathcal{B}\right) $. Then $\mathcal{A}$ is free.
\end{lem}

\begin{pf}
Let $\underline{e}=\sum_{i=1}^{m}b_{i}\underline{a_{i}}=\sum_{i=1}^{k}d_{i} 
\underline{c_{i}}$ with $\underline{a_{1}},\ldots ,\underline{a_{m}}, 
\underline{c_{1}},\ldots ,\underline{c_{k}}\in \mathcal{A}$ and $%
b_{1},\ldots ,b_{m},d_{1},\ldots ,d_{k}\in \mathcal{B}\backslash \left\{
0\right\} $. Note that $d_{i}=\left\langle \underline{e},\underline{c_{i}}
\right\rangle $ for $i=1,\ldots ,k$. Now if $\underline{c_{j}}\not\in
\left\{ \underline{a_{1}},\ldots ,\underline{a_{m}}\right\} $ for some $j\in
\left\{ 1,\ldots ,k\right\} $ then $d_{j}=\left\langle \underline{c_{j}}, 
\underline{e}\right\rangle =\left\langle \underline{c_{j}}
,\sum_{i=1}^{m}b_{i}\underline{a_{i}}\right\rangle =0$ which is a
contradiction. Hence $\left\{ \underline{c_{1}},\ldots ,\underline{c_{k}}
\right\} \subseteq \left\{ \underline{a_{1}},\ldots ,\underline{a_{m}}
\right\} $. The reverse inclusion is obtained by symmetry. Then for all $%
i=1,\ldots ,m$ there exists $j\in \left\{ 1,\ldots ,m\right\} $ such that $%
b_{i}=\left\langle \underline{e},\underline{a_{i}}\right\rangle
=\left\langle \underline{e},\underline{c_{j}}\right\rangle =d_{j}$,
concluding this proof.\qed
\end{pf}

\qquad We thus can set:

\begin{defn}
A subset $\mathcal{A}$ of $\mathcal{L}_{n}\left( \mathcal{B}\right) $ is an
orthonormal basis of $\mathcal{L}_{n}\left( \mathcal{B}\right) $ when it is
an orthonormal generating subset of $\mathcal{L}_{n}\left( \mathcal{B}
\right) $.
\end{defn}

\qquad An orthonormal basis is thus a generating set which, by Lemma (\ref%
{OrthoFree}), is also free, so it is basis, so that our vocabulary is
consistent.

\qquad There always exist orthonormal bases of $\mathcal{L}_{n}\left( 
\mathcal{B}\right) $ and we now give some examples. First, the \emph{\
canonical basis }or \emph{standard basis }of $\mathcal{L}_{n}\left( \mathcal{%
\ B}\right) $ is defined as the basis $\left( \underline{\delta _{i}}\right)
_{i=1,\ldots ,n}$ with $\underline{\delta _{1}}=\left( 1,0,\ldots ,0\right) $
, $\underline{\delta _{2}}=\left( 0,1,0,\ldots ,0\right) $, \ldots , $%
\underline{\delta _{n}}=\left( 0,\ldots ,0,1\right) $. More generally, we
have:

\begin{exmp}
\label{CyclicBasis}Let $\underline{a}=\left( a_{1},\ldots ,a_{n}\right) $ be
a stochastic vector. Let 
\begin{equation*}
\underline{e_{i}}=\left( a_{i},a_{i+1},\ldots ,a_{n},a_{1},\ldots
,a_{i-1}\right)
\end{equation*}
for all $i\in \left\{ 1,\ldots ,n\right\} $. Then by construction, $\left( 
\underline{e_{i}}\right) _{i=1,\ldots ,n}$ is an orthonormal subset of $%
\mathcal{L}_{n}\left( \mathcal{B}\right) $. Moreover 
\begin{eqnarray*}
\underline{\delta _{1}} &=&a_{1}\underline{e_{1}}+a_{2}\underline{e_{2}}
+\ldots +a_{n}\underline{e_{n}} \\
\underline{\delta _{2}} &=&a_{2}\underline{e_{1}}+a_{3}\underline{e_{2}}
+\ldots +a_{n}\underline{e_{n-1}}+a_{1}\underline{e_{n}} \\
&&\vdots \\
\underline{\delta _{n}} &=&a_{n}\underline{e_{1}}+a_{1}\underline{e_{2}}
+\ldots +a_{n-1}\underline{e_{n}}
\end{eqnarray*}
so $\left( \underline{e_{i}}\right) _{i=1,\ldots ,n}$ is a generating set
and thus an orthonormal basis for $\mathcal{L}_{n}\left( \mathcal{B}\right) $
.
\end{exmp}

\bigskip \qquad Let us observe that in general, linear independence in $%
\mathcal{L}_{n}\left( \mathcal{B}\right) $ is not an easy concept. We
propose in this paper to use orthogonality as a substitute. Indeed, if $%
\left\{ v_{1},\ldots ,v_{k}\right\} $ is a generating subset of $\mathcal{L}%
_{n}\left( \mathcal{B}\right) $ made of pairwise orthogonal, nonzero
vectors, then it is a minimal generating set, in the sense that any strict
subset is not generating (since, say, $v_{i}$ is not a linear combination of
the vectors in $\left\{ v_{1},\ldots ,v_{k}\right\} \backslash \left\{
v_{i}\right\} $ as all such combinations are orthogonal to $v_{i}$, the
inner product is definite yet $v_{i}\not=0$). However, orthogonality still
allows for some pathologies. For instance, assume there exists $a\in 
\mathcal{B}$ such that $a$ is neither $0$ or $1$. Then $\left( a,0\right) $, 
$\left( a^{c},0\right) $ and $\left( 0,1\right) $ are three nonzero
orthogonal vectors generating $\mathcal{L}_{2}\left( \mathcal{B}\right) $.
It is a minimal generating set, yet its cardinality is not minimal among all
generating families (since the canonical basis of $\mathcal{L}_{2}\left( 
\mathcal{B}\right) $ has cardinal $2$). If $\mathcal{B}$ is large enough, we
can even build on the same model infinite orthogonal generating families of
nonzero vectors, which are therefore minimal! We shall prove in the next
section that these pathologies are avoided when one restricts one's
attention to orthonormal bases. We shall also see that the concept of a
basis, i.e. a free generating subset, is in fact identical to the concept of
an orthonormal basis.

\bigskip \qquad The natural maps for our structure are:

\begin{defn}
A map $T:\mathcal{L}_{n}\left( \mathcal{B}\right) \longrightarrow \mathcal{L}
_{m}\left( \mathcal{B}\right) $ is linear when for all $a\in \mathcal{B}, 
\underline{b},\underline{c}\in \mathcal{L}_{n}\left( \mathcal{B}\right) $ we
have $T(a\underline{b}+\underline{c})=aT(\underline{b})+T(\underline{c})$.
\end{defn}

\bigskip \qquad As usual, $T(0)=0$ when $T$ is linear. When $T$ is linear
from $\mathcal{L}_{n}\left( \mathcal{B}\right) $ into $\mathcal{L}_{n}\left( 
\mathcal{B}\right) $, we call $T$ an operator on $\mathcal{L}_{n}\left( 
\mathcal{B}\right) $. An operator $T$ on $\mathcal{L}_{n}\left( \mathcal{B}
\right) $ is invertible when there exists an operator $S$ such that $S\circ
T=T\circ S=I$ where $I:x\in \mathcal{L}_{n}\left( \mathcal{B}\right) \mapsto
x$ is the identity operator. In the usual way, one can check that $T$ is an
invertible operator if and only if $T$ is a linear bijection, and the
inverse is a unique operator and is denoted by $T^{-1}$.

\bigskip \qquad We shall denote by $\mathcal{B}^{n}$ the Boolean algebra
product of $\mathcal{B}$ with itself $n$ times. Of course, the elements of $%
\mathcal{B}^{n}$ are the same as the elements of $\mathcal{L}_{n}\left( 
\mathcal{B}\right) $, but the algebraic structures are different.

\begin{lem}
\label{AutoAuto}If $T$ is an invertible operator on $\mathcal{L}_{n}\left( 
\mathcal{B}\right) $ then $T$ is a Boolean algebra automorphism on $\mathcal{%
\ B}^{n}$.
\end{lem}

\begin{pf}
Note that the supremum operation $\vee $ on $\mathcal{B}^{n}$ agrees with
the addition on $\mathcal{L}_{n}\left( \mathcal{B}\right) $ by definition.
So for any operator $L$ on $\mathcal{L}_{n}\left( \mathcal{B}\right) $ we
have $L\left( \underline{a}\vee \underline{b}\right) =L(\underline{a})\vee
L( \underline{b})$ and $L$ preserves the order $\leq $ on $\mathcal{B}^{n}$.
Hence, $T$ and $T^{-1}$ both preserve the order. Consequently, $\underline{a}
\leq \underline{b}$ if and only if $T(\underline{a})\leq T(\underline{b})$.
Hence $T$ is a lattice morphism, i.e. it also preserves the infimum. Also
note that this implies that $T(1,\ldots ,1)=\left( 1,\ldots ,1\right) $ --
since $\left( 1,\ldots ,1\right) $ is the largest element of $\mathcal{B}
^{n} $, we deduce that $T$ preserves the complement operation as well. This
concludes the proof.\qed
\end{pf}

\bigskip \qquad The converse of Lemma (\ref{AutoAuto}) does not hold,
namely: if $T:\mathcal{B}^{n}\longrightarrow \mathcal{B}^{n}$ is a Boolean
algebra automorphism then $T:\mathcal{L}_{n}\left( \mathcal{B}\right)
\longrightarrow \mathcal{L}_{n}\left( \mathcal{B}\right) $ need not be
linear. For example, let $\mathcal{B}=\left\{ 0,1,\omega ,\omega
^{c}\right\} $ and consider the Boolean algebra $\mathcal{B}^{2}$. Define
the automorphism $S$ on $\mathcal{B}$ by $S(\omega )=\omega ^{c}$ (so that $%
S(0)=0$, $S(1)=1$ and $S(\omega ^{c})=\omega $). Then $T=S\times S$ is an
automorphism of $\mathcal{B}^{2}$. Yet, seen as a map on $\mathcal{L}
_{2}\left( \mathcal{B}\right) $ we have 
\begin{equation*}
T\left( \omega \left( 1,0\right) \right) =T(\omega ,0)=\left( \omega
^{c},0\right)
\end{equation*}
and yet 
\begin{equation*}
\omega T\left( 1,0\right) =\left( \omega ,0\right)
\end{equation*}
and thus $T$ is not linear.

\bigskip \qquad We now show that if $\mathcal{B}$ is a finite Boolean
algebra, then any orthonormal basis for $\mathcal{L}_{n}\left( \mathcal{B}
\right) $ has cardinality $n$. Indeed, let $\left\{ \underline{e_{1}},\ldots
,\underline{e_{m}}\right\} $ be an orthonormal basis for $\mathcal{L}
_{n}\left( \mathcal{B}\right) $. Define $T:\mathcal{L}_{n}\left( \mathcal{B}
\right) \longrightarrow \mathcal{L}_{m}\left( \mathcal{B}\right) $ by 
\begin{equation*}
T\left( \underline{a}\right) =\left( \left\langle \underline{a},\underline{
e_{1}}\right\rangle ,\ldots ,\left\langle \underline{a},\underline{e_{m}}
\right\rangle \right) \text{.}
\end{equation*}
Then $T$ is a bijection from $\mathcal{B}^{n}$ onto $\mathcal{B}^{m}$ by
definition of orthonormal basis. Hence $n=m$ since $\mathcal{B}$ is finite.
As previously mentioned, we shall show in the next section that this result
holds for any Boolean algebra $\mathcal{B}$. Also, notice that $T$ thus
defined is an invertible operator on $\mathcal{L}_{n}\left( \mathcal{B}
\right) $, hence a Boolean algebra automorphism of $\mathcal{B}^{n}$ by
Lemma\ (\ref{AutoAuto}).

\bigskip \qquad As in traditional linear algebra, the study of linear maps
is facilitated by introducing matrices. A \emph{Boolean matrix} $A$ is a $%
n\times m$ matrix with entries in $\mathcal{B}$. We then write $A=\left[
a_{ij}\right] $ with $a_{ij}\in \mathcal{B}$ for $i\in \left\{ 1,\ldots
,n\right\} $ and $j\in \left\{ 1,\ldots ,m\right\} $. If $A$ is an $n\times
m $ Boolean matrix and if $B$ is an $m\times k$ Boolean matrix, then we
define the product $AB$ as the $n\times k$ matrix whose $\left( i,j\right) $
entry is given by $\vee _{p=1}^{m}a_{ip}b_{pj}$. In particular, we see
elements of $\mathcal{L}_{n}\left( \mathcal{B}\right) $ as $n\times 1$
matrices (i.e. column vectors). Boolean matrices, and a generalization to
distributive lattices have a considerable literature of investigation \cite%
{Blyth67,Cechlarova03,Giveon64,Jagannadham66,Luce52,Rutherford63,Rutherford63b,Rutherford64,Sindak75,Skornyakov86,Wedderburn34,Tan98,Tan01,Yoeli61}
.. These matrices provide useful tools in various fields such as switching
nets, automata theory and finite graph theory. Notice that permutation
matrices are a special case of (invertible) Boolean matrices.

\bigskip \qquad Our main motivation for studying Boolean matrices comes from
an analogy of a Markov chain \cite{Gudder08,Gudder08b,Stirzaker05}. Let $G$
be a finite directed graph whose vertices are labelled $1,2,\ldots ,n$ and
let $\mathcal{B}$ be a fixed Boolean algebra. We think of the vertices of $G$
as sites that a physical system can occupy. The edges of $G$ designate the
allowable transitions between sites. If there is an edge from vertex $i$ to
vertex $j$, we label it by an element $a_{ji}$ of $\mathcal{B}$. We think of 
$a_{ji}$ as the event, or proposition that the system evolves from site $i$
to site $j$ in one time-step. If there is no edge between $i$ and $j$ then
we set $a_{ji}=0$. The Boolean matrix $A=\left[ a_{ij}\right] $ is the
transition matrix in one-time-step for the physical system. The transition
matrix for $m$-time-steps is then naturally given by $A^{m}$.

\qquad Assuming that the system evolves from a site $i$ to some specific
site $j$ in one-time-step, we postulate that $a_{ji}a_{ki}=0$ for $j\not=k$
and $\vee _{j=1}^{n}a_{ji}=1$ for all $i=1,\ldots ,n$. Thus each column of $%
A $ is a stochastic vector. In the next section, we will refer to such
matrices as stochastic matrices. Suppose that $b_{i}$ is the event that the
system is in the site $i$ initially. We would then have that the vector $%
\underline{b}=\left( b_{1},\ldots ,b_{n}\right) $ is a stochastic vector and 
$A\underline{b}$ describes the system location after one-time-step. As we
shall see, $A\underline{b}$ is again a stochastic vector and in a natural
way, $\left( A\underline{b}\right) _{i}=\dbigvee\limits_{j=1}^{n}a_{ij}b_{j}$
is the event that the system is at site $i$ at one time-step. Thus, $m\in 
\mathbb{N}\mapsto A^{m}$ describes the dynamics of the system and this is
analogous to a traditional Markov chain. If in addition, we impose the
condition that for every site $i$ there is a specific site $j$ from which
the system evolved in one time-step, then we would have $a_{ij}a_{ik}=0$ and 
$\dbigvee\limits_{j=1}^{n}a_{ij}=1$. Such matrices are called unitary and
will be studied from Section 4 onward.

\qquad In general, if $G$ is a directed graph with $n$ vertices and $A$ is
an $n\times n$ stochastic matrix corresponding to the edges of $G$, we call $%
\left( G,A\right) $ a Boolean Markov chains. In section 6, we study the
powers of $A$ which are important for the description of the dynamics of $%
\left( G,A\right) $.

\section{The Dimension Theorem}

\qquad An orthonormal set is said to be \emph{stochastic} if all of its
elements are stochastic. In this section, we show that all orthonormal bases
of $\mathcal{L}_{n}\left( \mathcal{B}\right) $ have cardinality $n$.
Conversely, we show that any stochastic orthonormal set with cardinality $n$
is a basis for $\mathcal{L}_{n}\left( \mathcal{B}\right) $.

\qquad We shall use the following notations. Given a set $\mathcal{A=}
\left\{ \underline{a}_{1},\ldots ,\underline{a}_{m}\right\} $ of $m$
vectors, we use the notation $\underline{a}_{j}=\left( a_{1j},\ldots
,a_{nj}\right) $ with $a_{ij}\in \mathcal{B}$ ($i=1,\ldots ,n$ and $%
j=1,\ldots ,m$). Thus, we often think about a set $\left\{ \underline{a}
_{1},\ldots ,\underline{a}_{m}\right\} $ as a matrix $\left[ a_{ij}\right]
_{n\times m}$ whose columns are the elements of the set. By abuse of
notation, we denote this matrix by $\mathcal{A}$ again.

\bigskip \qquad We first establish that orthonormal bases possess a duality
property

\begin{thm}
\label{Duality}Let $\mathcal{A=}\left\{ \underline{a}_{1},\ldots ,\underline{%
a}_{m}\right\} $ be an orthonormal subset of $\mathcal{L}_{n}\left( \mathcal{%
B}\right) $. Then $\mathcal{A}$ is an orthonormal basis for $\mathcal{L}%
_{n}\left( \mathcal{B}\right) $ if and only if the set $\mathcal{A}^{\ast }$
of columns of $\left[ a_{ji}\right] _{m\times n}$ is an orthonormal subset
of $\mathcal{L}_{m}\left( \mathcal{B}\right) $.
\end{thm}

\begin{pf}
For all $j\in \left\{ 1,\ldots ,m\right\} $ we denote $\underline{a_{j}}
=\left( a_{1j},\ldots ,a_{nj}\right) $. Assume that $\mathcal{A}$ is an
orthonormal basis for $\mathcal{L}_{n}\left( \mathcal{B}\right) $. Then
there exists $b_{1},\ldots ,b_{m}\in \mathcal{B}$ such that $\underline{
\delta _{1}}=\sum_{j=1}^{m}b_{j}\underline{a_{j}}$. In particular, $%
0=\dbigvee\limits_{j=1}^{m}b_{j}a_{ij}$ for $i\not=1$ so $b_{j}a_{ij}=0$ for
all $i\in \left\{ 2,\ldots ,n\right\} $ and all $j\in \left\{ 1,\ldots
,m\right\} $. Hence, $b_{j}a_{1j}=b_{j}\left(
\dbigvee\limits_{i=1}^{n}a_{ij}\right) =b_{j}$ since $\vee
_{i=1}^{n}a_{ij}=1 $. Hence $b_{j}\leq a_{1j}$ for all $j\in \left\{
1,\ldots ,m\right\} $. On the other hand, $1=\dbigvee
\limits_{j=1}^{m}b_{j}a_{1j}$ and $a_{1j}$ and $a_{1k}$ are disjoint for $%
j\not=k$, so we must have $b_{j}a_{1j}=a_{1j}$ for all $j\in \left\{
1,\ldots ,m\right\} $. Consequently, $\dbigvee\limits_{j=1}^{m}a_{1j}=1$.
Moreover since $b_{j}a_{ij}=0$ for $i\not=1$, we conclude that $%
a_{1j}a_{ij}=0$ for $i\not=1$.

Replacing $\underline{\delta _{1}}$ by $\underline{\delta _{k}}$ for $k\in
\left\{ 1,\ldots ,n\right\} $ we see similarly that $\dbigvee
\limits_{j=1}^{m}a_{kj}=1$ and $a_{kj}a_{ij}=0$ for $i\not=k$ and for all $%
j\in \left\{ 1,\ldots ,m\right\} $. Hence, the set of columns of $\left[
a_{ji}\right] _{m\times n}$ is indeed an orthonormal subset of $\mathcal{L}
_{m}\left( \mathcal{B}\right) $.

Conversely, assume that $\mathcal{A}^{\ast }$ is an orthonormal subset of $%
\mathcal{L}_{m}\left( \mathcal{B}\right) $. This means by definition, and
using the same notations as before, that $\dbigvee\limits_{j=1}^{m}a_{ij}=1$
for all $i=1,\ldots ,n$ and $a_{kj}a_{ij}=0$ for all $i\not=k$ between $1$
and $n$ and $j=1,\ldots ,m$. It follows that 
\begin{equation}
\dbigvee\limits_{j=1}^{m}a_{kj}a_{ij}=\delta _{ik}\text{\ \ \ \ }%
(k,i=1,\ldots ,n)  \label{DualityEq1}
\end{equation}%
where $\delta _{ij}$ is $1\in \mathcal{B}$ if $i=j$ and $0\in \mathcal{B}$
otherwise. Now (\ref{DualityEq1}) is equivalent to 
\begin{equation*}
\underline{\delta _{k}}=\dbigvee\limits_{j=1}^{m}a_{kj}\underline{a_{j}}
\end{equation*}%
for $k=1,\ldots ,n$ and thus $\left\{ \underline{a_{1}},\ldots ,\underline{%
a_{m}}\right\} $ generates $\mathcal{L}_{n}\left( \mathcal{B}\right) $ and,
since it is an orthonormal set by assumption, it is an orthonormal basis of $%
\mathcal{L}_{n}\left( \mathcal{B}\right) $.\qed
\end{pf}

\begin{cor}
\label{StochasticBasis}An orthonormal basis is stochastic.
\end{cor}

\begin{cor}
\label{NSize}If $\left\{ \underline{a_{1}},\ldots ,\underline{ a_{n}}
\right\} $ is a stochastic orthonormal subset of $\mathcal{L} _{n}\left( 
\mathcal{B}\right) $ then it is a basis.
\end{cor}

\begin{pf}
Let $a=\left( \dbigvee\limits_{j=1}^{n}a_{1j}\right) ^{c}$ and assume $%
a\not=0$. By Stone's Theorem, there exists a set $\Omega ,$ a Boolean
algebra of subsets of $\Omega $ and a Boolean algebra isomorphism $\mathcal{%
B }\longrightarrow \mathcal{B}_{\Omega }$. We identify $\mathcal{B}$ and $%
\mathcal{B}_{\Omega }$ in this proof and thus regard the elements of $%
\mathcal{B}$ as subsets of $\Omega $, with $0$ identified with $\emptyset $
and $1$ with $\Omega $.

Let $\omega \in a$. Then $\omega \not\in a_{1j}$ for $j=1,\ldots ,n$.\ Since 
$\mathcal{A}$ is stochastic and orthonormal, we must have that $\omega \in
a_{i_{1}1}$, $\omega \in a_{i_{2}2}$, \ldots , $\omega \in a_{i_{n-1}n-1}$
for some $i_{1},\ldots ,i_{n-1}$ with $i_{r}\not=1$ and $i_{r}\not=i_{s}$
for $r,s=1,\ldots ,n-1$. Now, suppose $\omega \in a_{kn}$ for some $k\in
\left\{ 1,\ldots ,n\right\} $. Then $k\not=1$ (since $\omega \in a$) and $%
k\not=i_{r}$ for $r=1,\ldots ,n-1$ (orthogonality). But this is a
contradiction since this precludes $n$ values for $k$ which can only take $n$
values. Hence $\omega \not\in a_{kn}$ for all $k\in \left\{ 1,\ldots
,n\right\} $. This contradicts, in turn, that $\underline{a_{n}}$ is a unit
vector, i.e. form a partition of $\Omega $. Hence, $a=0$.

The same reasoning applies to show that $\dbigvee\limits_{j=1}^{n}a_{kj}=1$
for all $k\in \left\{ 1,\ldots ,n\right\} $. Hence $\mathcal{A}^{\ast }$ is
an orthonormal subset of $\mathcal{L}_{n}\left( \mathcal{B}\right) $ and
thus by Theorem (\ref{Duality}), $\mathcal{A}$ is an orthonormal basis for $%
\mathcal{L}_{n}\left( \mathcal{B}\right) $.\qed
\end{pf}

\bigskip \qquad By symmetry, we can restate Theorem (\ref{Duality}) by
stating that $\mathcal{A}$ is an orthonormal basis for $\mathcal{L}
_{n}\left( \mathcal{B}\right) $ if and only if $\mathcal{A}^{\ast }$ is an
orthonormal basis for $\mathcal{L}_{m}\left( \mathcal{B}\right) $. We call $%
\mathcal{A}^{\ast }$ the dual basis for $\mathcal{A}$. For example, if $%
a_{1},a_{2},a_{3}\in \mathcal{B}$ with $a_{1}\vee a_{2}\vee a_{3}=1$ and $%
a_{i}a_{j}=0$ for $i\not=j$ in $\left\{ 1,2,3\right\} $, then the columns of
the following matrix: 
\begin{equation*}
\mathcal{A=}\left[ 
\begin{array}{ccc}
a_{1} & a_{3} & a_{2} \\ 
a_{2} & 0 & a_{2}^{c} \\ 
a_{3} & a_{3}^{c} & 0%
\end{array}
\right]
\end{equation*}
form an orthonormal basis for $\mathcal{L}_{3}\left( \mathcal{B}\right) $.
The rows form the corresponding dual basis. Notice that $\mathcal{A}$ need
not be symmetric. Such a matrix $\mathcal{A}$ is what we shall call a
unitary matrix in section 4.

\bigskip \qquad We now establish a core result concerning the construction
of stochastic vectors.

\begin{thm}
\label{Descent}Let $n>1$. Let $\underline{a}=\left( a_{1},\ldots
,a_{n}\right) $ and $\underline{b}=\left( b_{1},\ldots ,b_{n}\right) $ be
two stochastic vectors in $\mathcal{L}_{n}\left( \mathcal{B}\right) $. Then $%
\underline{a}\perp \underline{b}$ if and only if there exists a stochastic
vector $\underline{c}=\left( c_{1},\ldots ,c_{n-1}\right) $ in $\mathcal{L}
_{n-1}\left( \mathcal{B}\right) $ such that $b_{i}=c_{i}a_{i}^{c}$ for $%
i=1,\ldots ,n-1$. If $\underline{a}\perp \underline{b}$ then we can always
choose $\underline{c}$ with $c_{i}=b_{n}a_{i}\vee b_{i}$ for $i=1,\ldots
,n-1 $.
\end{thm}

\begin{pf}
Suppose that $\underline{a}\perp \underline{b}$. Let $i\in \left\{ 1,\ldots
,n-1\right\} $. We set $c_{i}=b_{n}a_{i}\vee b_{i}$. Since $\underline{a}
\perp \underline{b}$, we have $b_{i}\leq a_{i}^{c}.$ Hence 
\begin{equation*}
c_{i}a_{i}^{c}=\left( b_{n}a_{i}\vee b_{i}\right)
a_{i}^{c}=b_{i}a_{i}^{c}=b_{i}\text{.}
\end{equation*}
Now, since $\underline{a}$ and $\underline{b}$ are stochastic vectors, we
conclude that for all $j\in \left\{ 1,\ldots ,n\right\} $ and $j\not=i$ we
have 
\begin{eqnarray*}
c_{i}c_{j} &=&\left( b_{n}a_{i}\vee b_{i}\right) \left( b_{n}a_{j}\vee
b_{j}\right) \\
&=&b_{n}a_{i}a_{j}\vee b_{n}b_{j}a_{i}\vee b_{i}b_{n}a_{j}\vee b_{i}b_{j}=0 
\text{.}
\end{eqnarray*}
Finally, we have 
\begin{eqnarray*}
\dbigvee\limits_{i=1}^{n-1}c_{i} &=&\dbigvee\limits_{i=1}^{n-1}\left(
b_{n}a_{i}\vee b_{i}\right) =\left(
b_{n}\dbigvee\limits_{i=1}^{n-1}a_{i}\right) \vee
\dbigvee\limits_{i=1}^{n-1}b_{i} \\
&=&b_{n}a_{n}^{c}\vee b_{n}^{c}=b_{n}\vee b_{n}^{c}=1\text{.}
\end{eqnarray*}
We conclude that $\underline{c}=\left( c_{1},\ldots ,c_{n-1}\right) $ is a
stochastic vector, and it obviously has the desired property.

Conversely, suppose that there exists a stochastic vector $\underline{c}$ in 
$\mathcal{L}_{n-1}\left( \mathcal{B}\right) $ such that $%
b_{i}=c_{i}a_{i}^{c} $ for $i=1,\ldots ,n-1$. Then by construction $%
a_{i}b_{i}=0$ for $i=1,\ldots ,n-1$. Moreover 
\begin{eqnarray*}
a_{n}b_{n} &=&a_{n}\left( \dbigvee\limits_{i=1}^{n-1}b_{i}\right)
^{c}=a_{n}\left( \dbigvee\limits_{i=1}^{n-1}c_{i}a_{i}^{c}\right) ^{c} \\
&=&a_{n}\dbigwedge\limits_{i=1}^{n-1}\left( a_{i}\vee c_{i}^{c}\right)
=a_{n}\dbigwedge\limits_{i=1}^{n-1}c_{i}^{c}=a_{n}\left(
\dbigvee\limits_{i=1}^{n-1}c_{i}\right) ^{c}=0\text{.}
\end{eqnarray*}
It follows that $\underline{a}\perp \underline{b}$.\qed
\end{pf}

\bigskip \qquad We can now show:

\begin{lem}
\label{BasisUpBound}If $\mathcal{A}=\left\{ \underline{a_{1}},\ldots , 
\underline{a_{m}}\right\} $ is a stochastic orthonormal set in $\mathcal{L}
_{n}\left( \mathcal{B}\right) $ then $m\leq n$.
\end{lem}

\begin{pf}
We proceed by induction on $n\in \mathbb{N}$. For $n=1$ the only orthonormal
set is $\left\{ 1\right\} $ so the result holds trivially. Now we assume the
results holds for some $n\in \mathbb{N}$. Let $\mathcal{A=}\left[ a_{ij}%
\right] _{(n+1)\times m}$ be a stochastic orthonormal set in $\mathcal{L}%
_{n+1}\left( \mathcal{B}\right) $. By Theorem (\ref{Descent}), for each $%
j=2,\ldots ,m$ there exists a stochastic vector $\underline{c_{j}}=\left(
c_{1j},\ldots ,c_{nj}\right) $ in $\mathcal{L}_{n}\left( \mathcal{B}\right) $
such that $a_{ij}=c_{ij}a_{i1}^{c}$ for all $i=1,\ldots ,n$ and $j=2,\ldots
,m$. Let $j,k\in \left\{ 2,\ldots ,m\right\} $ with $j\not=k$ and $i\in
\left\{ 1,\ldots ,n\right\} $. Recall from Theorem (\ref{Descent}) that $%
c_{ij}=a_{i1}a_{nj}\vee a_{ij}$, and since $\mathcal{A}$ is orthonormal 
\begin{eqnarray*}
c_{ij}c_{ik} &=&\left( a_{i1}a_{n+1,j}\vee a_{ij}\right) \left(
a_{i1}a_{n+1,k}\vee a_{ik}\right) \\
&=&a_{i1}a_{n+1,j}a_{n+1,k}\vee a_{i1}a_{n+1,j}a_{ik}\vee
a_{i1}a_{ij}a_{n+1,k}\vee a_{ij}a_{ik} \\
&=&0\text{.}
\end{eqnarray*}

Hence $\left\{ \underline{c_{2}},\ldots ,\underline{c_{m}}\right\} $ is a
stochastic orthonormal set in $\mathcal{L}_{n}\left( \mathcal{B}\right) $.
By our induction hypothesis, $m-1\leq n$ and thus $m\leq n+1$, which
completes our proof by induction.\qed
\end{pf}

\qquad The main result of this section is:

\begin{thm}
\label{Dimension}If $\mathcal{A}$ is an orthonormal basis for $\mathcal{L}
_{n}\left( \mathcal{B}\right) $ then the cardinality of $\mathcal{A}$ is $n$.
\end{thm}

\begin{pf}
We proceed by induction on $n$. The result is trivial for $n=1$. Assume that
for some $n\in \mathbb{N}$, if $\mathcal{A}_{0}$ is an orthonormal basis for 
$\mathcal{L}_{k}\left( \mathcal{B}\right) $ with $k\leq n$ then $\mathcal{A}%
_{0}$ contains exactly $k$ vectors. Let $\mathcal{A}$ be an orthonormal
basis of $\mathcal{L}_{n+1}\left( \mathcal{B}\right) $. By Corollary (\ref%
{StochasticBasis}), $\mathcal{A}$ is stochastic. Applying Lemma (\ref%
{BasisUpBound}), we deduce that the cardinality $m$ of $\mathcal{A}$
satisfies $m\leq n+1$. Assume that $m<n+1$. By Theorem (\ref{Duality}), $%
\mathcal{A}^{\ast }$ is an orthonormal basis for $\mathcal{L}_{m}\left( 
\mathcal{B}\right) $ since $\mathcal{A}=\left( \mathcal{A}^{\ast }\right)
^{\ast }$ is an orthonormal subset of $\mathcal{L}_{n+1}\left( \mathcal{B}%
\right) $. Since $m\leq n$, we conclude by our induction hypothesis that the
cardinality of $\mathcal{A}^{\ast }$ is $m$. But by construction, the
cardinality of $\mathcal{A}^{\ast }$ is $n+1$, which is a contradiction.
Hence $m=n+1$ which completes our proof by induction.\qed
\end{pf}

Combining Theorem (\ref{Dimension}) and Corollary (\ref{NSize}) we obtain
the following result:

\begin{cor}
\label{DimCor}A stochastic orthonormal set $\mathcal{A}$ is a basis for $%
\mathcal{L}_{n}\left( \mathcal{B}\right) $ if and only if the cardinality of 
$\mathcal{A}$ is $n$.
\end{cor}

\qquad To be fully satisfactory, we shall now check that the orthonormal
families of $\mathcal{L}_{n}\left( \mathcal{B}\right) $ of cardinality $n$
are in fact basis. We shall use the following:

\begin{lem}
\label{Redux1}If $\underline{a}=\left( a_{1},\ldots ,a_{n}\right) \in 
\mathcal{L}_{n}\left( \mathcal{B}\right) $ is a unit vector, then there
exists a stochastic vector $\underline{b}=\left( b_{1},\ldots ,b_{n}\right) $
with $b_{i}\leq a_{i}$ for all $i=1,\ldots ,n$.
\end{lem}

\begin{pf}
For $i=1,\ldots ,n$ we set $b_{i}=a_{i}\left( a_{1}^{c}a_{2}^{c}\ldots
a_{i-1}^{c}\right) \leq a_{i}$. Then $b_{i}b_{j}=0$ for $i,j=1,\ldots ,n$
and $i\not=j$, and $\dbigvee\limits_{i=1}^{n}b_{i}=\dbigvee
\limits_{i=1}^{n}a_{i}=1$ so $\underline{b}$ is a stochastic vector.\qed
\end{pf}

\qquad Now, we can state:

\begin{cor}
\label{DimCor2}An orthonormal set of $\mathcal{L}_{n}\left( \mathcal{B}
\right) $ is a basis for $\mathcal{L}_{n}\left( \mathcal{B}\right) $ if and
only if it has cardinality $n$.
\end{cor}

\begin{pf}
Let $\mathcal{A=}\left\{ \underline{a_{1}},\ldots ,\underline{a_{n}}\right\} 
$ be an orthonormal set. Using Lemma (\ref{Redux1}), there exists a set of
stochastic vectors $\underline{b_{1}},\ldots ,\underline{b_{n}}$ such that $%
b_{ij}\leq a_{ij}$. Therefore, $\left\{ \underline{b_{1}},\ldots ,\underline{
b_{n}}\right\} $ is a stochastic orthogonal set of size $n$ and thus it is a
basis for $\mathcal{L}_{n}\left( \mathcal{B}\right) $ by Corollary (\ref%
{DimCor}). Now, let $i,j,k,l=1,\ldots ,n$ with $i>j$. Let $\underline{v}=$ $%
a_{ik}a_{jk}\underline{\delta _{i}}$. Then, using the construction of Lemma
( \ref{Redux1}), we have 
\begin{equation*}
a_{ik}a_{jk}b_{il}=a_{ik}a_{jk}a_{il}a_{1l}^{c}\ldots a_{i-1,l}^{c}=0
\end{equation*}
since either $l=k$ and then $a_{ik}a_{jk}b_{il}\leq a_{jk}a_{jk}^{c}=0$
since $i>j$, or $l\not=k$ and $a_{ik}a_{il}=0$ since $\mathcal{A}$ is
orthogonal. Hence the vector $\underline{v}$ is orthogonal to $\underline{
b_{1}},\ldots ,\underline{b_{n}}$, thus $\underline{v}=0$. Hence, $\mathcal{%
A }$ is stochastic.\ By Corollary (\ref{DimCor}), it is an orthonormal basis
of $\mathcal{L}_{n}\left( \mathcal{B}\right) $.

The converse is Theorem (\ref{Dimension}).\qed
\end{pf}

\qquad In view of Corollary (\ref{DimCor2}), we call $n$ the \emph{dimension}
of the Boolean linear space $\mathcal{L}_{n}\left( \mathcal{B}\right) $. We
now consider the following question: can any stochastic orthonormal subset $%
\mathcal{A}$ of $\mathcal{L}_{n}\left( \mathcal{B}\right) $ be extended to
an orthonormal basis for $\mathcal{L}_{n}\left( \mathcal{B}\right) $? By
Lemma (\ref{BasisUpBound}), $\mathcal{A}$ can not have more than $n$
vectors. Of course, if the cardinality of $\mathcal{A}$ is $n$ then it is
already a basis by Corollary (\ref{NSize}). Moreover, Example (\ref%
{CyclicBasis}) shows that if $\mathcal{A}$ is reduced to a unique stochastic
vector, then there is an orthonormal basis of $\mathcal{L}_{n}\left( 
\mathcal{B}\right) $ containing $\mathcal{A}$ so the answer is affirmative.
We shall now prove that the answer is affirmative in general.

\bigskip \qquad We shall use the following concept:

\begin{defn}
A subset $\mathcal{M}\subseteq \mathcal{L}_{n}\left( \mathcal{B}\right) $ is
a subspace if it is generated by an orthonormal set $\mathcal{A}=\left\{ 
\underline{a_{1}},\ldots ,\underline{a_{m}}\right\} $, i.e. 
\begin{equation*}
\mathcal{M}=\left\{ \sum_{i=1}^{m}b_{i}\underline{a_{i}}:b_{1},\ldots
,b_{m}\in \mathcal{B}\right\} \text{.}
\end{equation*}

Any orthonormal set $\mathcal{A}$ generating $\mathcal{M}$ is called an
orthonormal basis for $\mathcal{M}$.
\end{defn}

\qquad We emphasize that we do not require orthonormal bases of subspaces to
be stochastic. In fact, a subspace may not contain any stochastic
orthonormal basis: for example, if there exists $a\in \mathcal{B}$ such that 
$a\not\in \left\{ 0,1\right\} $ then the subset $E=\left\{ b\left(
1,a\right) :b\in \mathcal{B}\right\} $ is a subspace with basis $\left(
1,a\right) $. Since any orthonormal set of two vectors generates $\mathcal{L}
_{2}\left( \mathcal{B}\right) \not=E$, any orthonormal basis for $E$ is
necessarily reduced to one vector. If this vector is stochastic, then it is
of the form $\left( b,b^{c}\right) $ for some $b\in \mathcal{B}$. It is then
easy to check that $(1,a)$ can not be of the form $\left( cb,cb^{c}\right) $
and thus $E$ has no stochastic vector basis. Thus, we will sometimes use:

\begin{defn}
A subspace with a stochastic orthonormal basis is called a stochastic
subspace.
\end{defn}

\qquad Linear maps generalize trivially to linear maps between two
subspaces. Of special interest to us will be:

\begin{defn}
A linear map $T:\mathcal{M}\longrightarrow \mathcal{N}$ between two
subspaces $\mathcal{M}$ and $\mathcal{N}$ of, respectively, $\mathcal{L}
_{n}\left( \mathcal{B}\right) $ and $\mathcal{L}_{m}\left( \mathcal{B}
\right) $, is called an isometry when for all $\underline{a},\underline{b}
\in \mathcal{M}$ we have $\left\langle T(\underline{a}),T(\underline{b}
)\right\rangle =\left\langle \underline{a},\underline{b}\right\rangle $.
\end{defn}

\begin{lem}
\label{IsometryChar}Let $\mathcal{M}\subseteq \mathcal{L}_{n}\left( \mathcal{%
\ B}\right) $ and $\mathcal{N}\subseteq \mathcal{L}_{m}\left( \mathcal{B}
\right) $ be two subspaces. Let $T:\mathcal{M}\longrightarrow \mathcal{N}$
be a linear map. The following are equivalent:

\begin{enumerate}
\item $T$ is an isometry,

\item There exists an orthonormal basis $\mathcal{A}=\left\{ \underline{%
e_{1} },\ldots ,\underline{e_{k}}\right\} $ of $\mathcal{M}$ such that $%
\left\{ T \underline{e_{i}}:i=1,\ldots ,k\right\} $ is an orthonormal set of 
$\mathcal{\ N}$,

\item For every orthonormal set $\mathcal{A=}\left\{ \underline{e_{1}}
,\ldots ,\underline{e_{k}}\right\} $ of $\mathcal{M}$, the set $\left\{ T 
\underline{e_{1}},\ldots T\underline{e_{k}}\right\} $ is an orthonormal set
of $\mathcal{N}$.
\end{enumerate}

Moreover, if $T$ is an isometry, then it is injective.
\end{lem}

\begin{pf}
We start by proving that (2) implies (1). Let $\mathcal{A=}\left\{ 
\underline{e_{1}},\ldots ,\underline{e_{k}}\right\} $ be an orthonormal
basis of $\mathcal{M}$ such that $\left\{ T\underline{e_{1}},\ldots ,T%
\underline{e_{k}}\right\} $ is orthonormal. Let $\underline{a},\underline{b}%
\in \mathcal{M}$. We can write $\underline{a}=\sum_{i=1}^{k}a_{i}\underline{%
e_{i}}$ and $\underline{b}=\sum_{i=1}^{k}b_{i}\underline{e_{i}}$ with $%
a_{i},b_{i}\in \mathcal{B}$ ($i=1,\ldots ,k$). Then 
\begin{eqnarray*}
\left\langle T\underline{a},T\underline{b}\right\rangle
&=&\dbigvee\limits_{i,j=1}^{k}\left\langle a_{i}T\underline{e_{i}},b_{j}T%
\underline{e_{j}}\right\rangle
=\dbigvee\limits_{i,j=1}^{k}a_{i}b_{j}\left\langle T\underline{e_{i}},T%
\underline{e_{j}}\right\rangle \\
&=&\dbigvee\limits_{i=1}^{k}a_{i}b_{i}=\left\langle \underline{a},\underline{%
b}\right\rangle \text{.}
\end{eqnarray*}%
Hence $T$ is an isometry.

Now, (1) implies (3) and (3) implies (2) are both trivial.

Assume now that $T$ is an isometry. Assume $T\underline{a}=T\underline{b}$.
Then, using the same notations as above, we have 
\begin{equation*}
a_{i}=\left\langle a,\underline{e_{i}}\right\rangle =\left\langle T 
\underline{a},T\underline{e_{i}}\right\rangle =\left\langle T\underline{b},T 
\underline{e_{i}}\right\rangle =\left\langle \underline{b},\underline{e_{i}}
\right\rangle =b_{i}
\end{equation*}
for all $i=1,\ldots ,k$. Hence $\underline{a}=\underline{b}$.\qed
\end{pf}

\begin{defn}
Let $\mathcal{M}$ and $\mathcal{N}$ be two subspaces of respectively $%
\mathcal{L}_{m}\left( \mathcal{B}\right) $ and $\mathcal{L}_{n}\left( 
\mathcal{B}\right) $. A surjective isometry $T:\mathcal{M}\longrightarrow 
\mathcal{N}$ is called an isomorphism, and then $\mathcal{M}$ and $\mathcal{%
N }$ are called isomorphic subspaces.
\end{defn}

\bigskip \qquad It is clear that the inverse of an isomorphism is an
isomorphism, and that the composition of two isomorphisms is again an
isomorphism. It follows that isomorphic is an equivalence relation. It is
also an important observation that isomorphisms map orthonormal bases to
orthonormal bases: if $\left\{ \underline{a_{1}},\ldots ,\underline{a_{n}}%
\right\} $ is an orthonormal basis for a subspace $\mathcal{M}$ and $T:%
\mathcal{M}\longrightarrow \mathcal{N}$ is an isomorphism then $\left\{ T%
\underline{a_{1}},\ldots ,T\underline{a_{n}}\right\} $ is an orthonormal set
since $T$ is an isometry (Lemma (\ref{IsometryChar})). Moreover, if $%
\underline{b}\in \mathcal{N}$ then there exists $\underline{c}\in \mathcal{M}
$ such that $T(\underline{c})=\underline{b}$. Since $\underline{c}%
=\sum_{i=1}^{n}c_{i}\underline{a_{i}}$ for some $c_{1},\ldots ,c_{n}\in 
\mathcal{B}$ we conclude that $\underline{b}=\sum_{i=1}^{n}\underline{c_{i}}%
T(\underline{a_{i}})$. Hence $\left\{ T\underline{a_{1}},\ldots ,T\underline{%
a_{n}}\right\} $ is an orthonormal generating subset of $\mathcal{N}$, hence
a basis of $\mathcal{N}$.

\begin{thm}
\label{SubIso}If $\mathcal{M}$ is a subspace then there exists an $m\in 
\mathbb{N}$ and an isomorphism $T:\mathcal{M}\longrightarrow \mathcal{L}
_{m}\left( \mathcal{B}\right) $. Moreover $T$ can be chosen to take
stochastic vectors to stochastic vectors, and if $\mathcal{M}$ is a
stochastic subspace then $T$ can be chosen so that $T$ and $T^{-1}$ map
stochastic vectors to stochastic vectors.
\end{thm}

\begin{pf}
Let $\left\{ \underline{e_{1}},\ldots ,\underline{e_{m}}\right\} $ be an
orthonormal basis for $\mathcal{M}$ and let us denote the canonical basis of 
$\mathcal{L}_{m}\left( \mathcal{B}\right) $ by $\left\{ \underline{\delta
_{1}} ,\ldots ,\underline{\delta _{m}}\right\} $. We define $T:\mathcal{M}
\longrightarrow \mathcal{L}_{m}\left( \mathcal{B}\right) $ by setting for
all $\underline{a}\in \mathcal{M}$: 
\begin{equation*}
T\underline{a}=\left( \left\langle \underline{a},e_{1}\right\rangle ,\ldots
,\left\langle \underline{a},e_{m}\right\rangle \right) \text{.}
\end{equation*}
Then $T$ is linear and $T\underline{e_{i}}=\underline{\delta _{i}}$ for $%
i\in \left\{ 1,\ldots ,m\right\} $. By Lemma (\ref{IsometryChar}), $T$ is an
isometry and $T$ is surjective by construction (if $\underline{b}\in 
\mathcal{L}_{m}\left( \mathcal{B}\right) $ then $\underline{b}=\left(
b_{1},\ldots ,b_{m}\right) $ then $T\left( \sum_{i=1}^{m}b_{i}\underline{
e_{i}}\right) =\underline{b}$). So $T$ is an isomorphism.

Moreover, $T$ preserves stochastic vectors. Indeed, let $\underline{a}$ be a
stochastic vector in $\mathcal{M}$. Let $n\in \mathbb{N}$ such that $%
\mathcal{M}$ is a subspace of $\mathcal{L}_{n}\left( \mathcal{B}\right) $.
Denote by $\left\{ \underline{\delta _{1}^{\prime }},\ldots ,\underline{
\delta _{n}^{^{\prime }}}\right\} $ the canonical orthonormal basis of $%
\mathcal{L}_{n}\left( \mathcal{B}\right) $. For $i=1,\ldots ,m$ then 
\begin{equation*}
\left\langle \underline{a},\underline{e_{i}}\right\rangle =\left\langle 
\underline{a},\sum_{r=1}^{n}\left\langle \underline{e_{i}},\underline{\delta
_{r}^{\prime }}\right\rangle \underline{\delta _{r}^{\prime }}\right\rangle
=\dbigvee\limits_{r=1}^{n}\left\langle \underline{e_{i}},\underline{\delta
_{r}^{\prime }}\right\rangle \left\langle \underline{a},\underline{\delta
_{r}^{\prime }}\right\rangle \text{.}
\end{equation*}
Hence, for $i\not=j$ and $i,j=1,\ldots ,m$ we have 
\begin{eqnarray*}
\left\langle \underline{a},\underline{e_{i}}\right\rangle \left\langle 
\underline{a},\underline{e_{j}}\right\rangle
&=&\dbigvee\limits_{r,s=1}^{n}\left\langle \underline{e_{i}},\underline{
\delta _{r}^{\prime }}\right\rangle \left\langle \underline{a},\underline{
\delta _{r}^{\prime }}\right\rangle \left\langle \underline{e_{j}}, 
\underline{\delta _{s}^{\prime }}\right\rangle \left\langle \underline{a}, 
\underline{\delta _{s}^{\prime }}\right\rangle \\
&=&\dbigvee\limits_{r=1}^{n}\left\langle \underline{e_{i}},\underline{\delta
_{r}^{\prime }}\right\rangle \left\langle \underline{a},\underline{\delta
_{r}^{\prime }}\right\rangle \left\langle \underline{e_{j}},\underline{
\delta _{r}^{\prime }}\right\rangle \text{ since }\underline{a}\text{ is
stochastic} \\
&\leq &\dbigvee\limits_{r=1}^{n}\left\langle \underline{e_{i}},\underline{
\delta _{r}^{\prime }}\right\rangle \left\langle \underline{e_{j}}, 
\underline{\delta _{r}^{\prime }}\right\rangle =\left\langle \underline{%
e_{i} },\underline{e_{j}}\right\rangle =0\text{.}
\end{eqnarray*}
Hence, by definition, $T\underline{a}$ is stochastic.

Now, it is easy to check that $T^{-1}\left( a_{1},\ldots ,a_{m}\right)
=\sum_{k=1}^{m}a_{k}\underline{e_{k}}$. Assume that $\mathcal{M}$ is
stochastic and that the basis $\left\{ \underline{e_{1}},\ldots ,\underline{%
e_{n}}\right\} $ is stochastic. If $\left( a_{1},\ldots ,a_{m}\right) \in 
\mathcal{L}_{m}\left( \mathcal{B}\right) $ is stochastic, then for $%
r,s=1,\ldots ,m$: 
\begin{eqnarray*}
\left\langle \sum_{k=1}^{m}a_{k}\underline{e_{k}},\underline{\delta
_{r}^{\prime }}\right\rangle \left\langle \sum_{k=1}^{m}a_{k}\underline{e_{k}%
},\underline{\delta _{s}^{\prime }}\right\rangle
&=&\dbigvee\limits_{k,l=1}^{m}a_{k}a_{l}\left\langle e_{k},\underline{\delta
_{r}^{\prime }}\right\rangle \left\langle e_{l},\underline{\delta
_{s}^{\prime }}\right\rangle \\
&=&\dbigvee\limits_{k=1}^{m}a_{k}\left\langle e_{k},\underline{\delta
_{r}^{\prime }}\right\rangle \left\langle e_{k},\underline{\delta
_{s}^{\prime }}\right\rangle \text{ as }\underline{a}\text{ is stochastic} \\
&=&a_{k}\delta _{r}^{s}
\end{eqnarray*}%
with $\delta _{r}^{s}$ is the Kronecker symbol. Note that we used that by
definition, an orthonormal basis of a subspace is stochastic. Hence $%
T^{-1}\left( a_{1},\ldots ,a_{m}\right) $ is a stochastic vector as well.
Hence $T^{-1}$ maps stochastic vectors to stochastic vectors.\qed
\end{pf}

\begin{cor}
Any two orthonormal bases of a subspace $\mathcal{M}$ have the same
cardinality.
\end{cor}

\begin{pf}
Let $\mathcal{A=}\left\{ \underline{a_{1}},\ldots ,\underline{a_{m}}\right\} 
$ and $\mathcal{B=}\left\{ \underline{b_{1}},\ldots ,\underline{b_{n}}
\right\} $ be two orthonormal bases of $\mathcal{M}$. By Theorem (\ref%
{SubIso}), there exists isomorphisms $T:\mathcal{M}\longrightarrow \mathcal{%
L }_{m}\left( \mathcal{B}\right) $ and $S:\mathcal{M}\longrightarrow 
\mathcal{L }_{n}\left( \mathcal{B}\right) $. Hence $T\circ S^{-1}:\mathcal{L}%
_{n}\left( \mathcal{B}\right) \longrightarrow \mathcal{L}_{m}\left( \mathcal{%
B}\right) $ is an isomorphism. In particular, it maps orthonormal basis to
orthonormal basis. Hence $n=m$ by Theorem (\ref{Dimension}).\qed
\end{pf}

\bigskip \qquad We call the common cardinality of all orthonormal bases for
a subspace $\mathcal{M}$ the \emph{dimension of }$\mathcal{M}$. It follows
from Theorem (\ref{SubIso}) that if $\mathcal{M}$ has dimension $m$, then $%
\mathcal{M}$ is isomorphic to $\mathcal{L}_{m}\left( \mathcal{B}\right) $. A
source of examples of subspaces is given by:

\begin{prop}
\label{OrthoSub}For any $\underline{a_{1}}\in \mathcal{L}_{n}\left( \mathcal{%
\ B}\right) $ we denote by $\underline{a_{1}}^{\perp }$ the set 
\begin{equation*}
\left\{ \underline{b}\in \mathcal{L}_{n}\left( \mathcal{B}\right)
:\left\langle \underline{a_{1}},\underline{b}\right\rangle =0\right\} \text{%
. }
\end{equation*}
If $\underline{a_{1}}$ is stochastic then $\underline{a_{1}}^{\perp }$ is a
stochastic subspace of $\mathcal{L}_{n}\left( \mathcal{B}\right) $ of
dimension $n-1$.
\end{prop}

\begin{pf}
Using Example (\ref{CyclicBasis}), we extend the stochastic vector $%
\underline{a_{1}}$ to an orthonormal basis $\left\{ \underline{a_{1}},\ldots
,\underline{a_{n}}\right\} $ of $\mathcal{L}_{n}\left( \mathcal{B}\right) $.
If $\underline{b}\perp \underline{a_{1}}$ then, writing $\underline{b}
=\sum_{i=1}^{n}b_{i}\underline{a_{i}}$ we see that $\left\langle \underline{%
b },\underline{a_{1}}\right\rangle =0$ if and only if $b_{1}=0$. Hence 
\begin{equation*}
\underline{a_{1}}^{\perp }=\left\{ \sum_{i=2}^{n}b_{i}\underline{a_{i}}
:b_{2},\ldots ,b_{n}\in \mathcal{B}\right\}
\end{equation*}
is the subspace generated by the stochastic orthonormal set $\left\{ 
\underline{a_{2}}\ldots ,\underline{a_{n}}\right\} $ of cardinality $n-1$. 
\qed
\end{pf}

\bigskip \qquad We are now ready to show:

\begin{thm}
\label{Incomplete}If $\mathcal{A}=\left\{ \underline{a_{1}},\ldots , 
\underline{a_{m}}\right\} $ is a stochastic orthonormal set in $\mathcal{L}
_{n}\left( \mathcal{B}\right) $ with $m<n$ then $\mathcal{A}$ can be
extended to an orthonormal basis for $\mathcal{L}_{n}\left( \mathcal{B}
\right) $.
\end{thm}

\begin{pf}
We proceed by induction on $n$. The result is trivial for $n=1$. Assume that
for some $n\in \mathbb{N}$, any stochastic orthonormal set of cardinality $%
m<n$ in $\mathcal{L}_{n}\left( \mathcal{B}\right) $ can be extended to a
basis for $\mathcal{L}_{n}\left( \mathcal{B}\right) $. Let $\mathcal{A=}%
\left\{ \underline{a_{1}},\ldots ,\underline{a_{m}}\right\} $ be a
stochastic orthonormal subset of $\mathcal{L}_{n+1}\left( \mathcal{B}\right) 
$ with $m<n+1$. By Proposition (\ref{OrthoSub}) and Theorem (\ref{SubIso})
there exist an isomorphism $T:\underline{a_{1}}^{\perp }\longrightarrow 
\mathcal{L}_{n}\left( \mathcal{B}\right) $ such that $T$ and $T^{-1}$
preserve stochastic vectors. Moreover, $\left\{ \underline{a_{2}},\ldots ,%
\underline{a_{m}}\right\} \subseteq \underline{a_{1}}^{\perp }$. Let $%
\underline{b_{i}}\in \mathcal{L}_{n}\left( \mathcal{B}\right) $ be given by $%
T\underline{a_{i}}=\underline{b_{i}}$ for $i=2,\ldots ,m$. It follows from
Lemma (\ref{IsometryChar}) that $\left\{ \underline{b_{2}},\ldots ,%
\underline{b_{m}}\right\} $ is an orthonormal set in $\mathcal{L}_{n}\left( 
\mathcal{B}\right) $ of cardinal $m-1<n$. By our induction hypothesis, there
exist stochastic vectors $\underline{b_{m+1}},\ldots ,\underline{b_{n+1}}$
such that $\left\{ \underline{b_{2}},\ldots ,\underline{b_{n+1}}\right\} $
is an orthonormal basis for $\mathcal{L}_{n}\left( \mathcal{B}\right) $. By
Theorem (\ref{SubIso}), $\left\{ T^{-1}\underline{b_{2}},\ldots ,T^{-1}%
\underline{b_{n+1}}\right\} $ is a stochastic orthonormal set in $\mathcal{L}%
_{n+1}\left( \mathcal{B}\right) $ which is a basis for $\underline{a_{1}}%
^{\perp }$. Since $\underline{a_{i}}=T^{-1}\underline{b_{i}}$ for $%
i=2,\ldots ,m$, we conclude by Corollary (\ref{NSize}) that $\left\{ 
\underline{a_{1}},T^{-1}\underline{b_{2}},\ldots ,T^{-1}\underline{b_{n+1}}%
\right\} $ is an orthonormal basis of $\mathcal{L}_{n+1}\left( \mathcal{B}%
\right) $ which extends $\mathcal{A}$.\qed
\end{pf}

\bigskip \qquad It follows from Theorem (\ref{Incomplete}) that if $\mathcal{%
\ M}$ is a stochastic subspace of $\mathcal{L}_{n}\left( \mathcal{B}\right) $
then 
\begin{equation*}
\mathcal{M}^{\perp }=\left\{ \underline{b}\in \mathcal{L}_{n}\left( \mathcal{%
\ B}\right) :\forall \underline{a}\in \mathcal{M}\ \ \ \ \underline{b}\perp 
\underline{a}\right\}
\end{equation*}
is also a stochastic subspace and $\mathcal{L}_{n}\left( \mathcal{B}\right)
= \mathcal{M}+\mathcal{M}^{\perp }$. One can now study projection operators
and the order structure of subspaces but we leave this for later work.

\section{Stochastic and Unitary Matrices}

\bigskip \qquad In the sequel, a matrix on $\mathcal{L}_{n}\left( \mathcal{B}
\right) $ will mean an $n\times n$ Boolean matrix, and a vector in $\mathcal{%
\ L}_{n}\left( \mathcal{B}\right) $ will mean a Boolean vector and will be
identified with a $n\times 1$ column vector. Moreover, if $A$ is a matrix
then we denote the $(i,j)^{\text{th}}$ entry by $(A)_{ij}$, or simply $%
\left( A\right) _{i}$ if $A$ is a column vector.

\qquad Let $A$ be a matrix on $\mathcal{L}_{n}\left( \mathcal{B}\right) $.
Then the map $\underline{x}\in \mathcal{L}_{n}\left( \mathcal{B}\right)
\mapsto A\underline{x}$ is linear and will be identified with $A$. Indeed,
for all $\underline{b},\underline{c}\in \mathcal{L}_{n}\left( \mathcal{B}
\right) $, $c\in \mathcal{B}$, and $i=1,\ldots ,n$ we have 
\begin{equation*}
\left( A\left( c\underline{b}\right) \right)
_{i}=\dbigvee\limits_{j=1}^{n}a_{ij}\left( c\underline{b}\right)
_{j}=\dbigvee\limits_{j=1}^{n}a_{ij}cb_{j}=c\dbigvee
\limits_{j=1}^{n}a_{ij}b_{j}=c\left( A\underline{b}\right) _{i}
\end{equation*}
and 
\begin{eqnarray*}
\left( A\left( \underline{b}+\underline{c}\right) \right) _{i}
&=&\dbigvee\limits_{j=1}^{n}a_{ij}\left( \underline{b}+\underline{c}\right)
_{j}=\dbigvee\limits_{j=1}^{n}a_{ij}\left( b_{j}\vee c_{j}\right) \\
=\left( \dbigvee\limits_{j=1}^{n}a_{ij}b_{j}\right) \vee \left(
\dbigvee\limits_{j=1}^{n}a_{ij}c_{j}\right) \\
&=&\left( A\underline{b}\right) _{i}\vee \left( A\underline{c}\right)
_{i}=\left( A\underline{b}+A\underline{c}\right) _{i}\text{.}
\end{eqnarray*}

\qquad Conversely, any operator $T$ on $\mathcal{L}_{n}\left( \mathcal{B}
\right) $ can be represented by a matrix on $\mathcal{L}_{n}\left( \mathcal{%
B }\right) $ with respect to the canonical basis. Indeed, define $%
a_{ij}=\left\langle T\underline{\delta _{j}},\underline{\delta _{i}}
\right\rangle $ for all $i,j=1,\ldots ,n$. Then $T\underline{\delta _{j}}
=\sum_{i=1}^{n}a_{ij}\underline{\delta _{i}}$. Defining the matrix $A_{T}= %
\left[ a_{ij}\right] _{n\times n}$ we have 
\begin{equation*}
\left( A_{T}\underline{\delta _{i}}\right)
_{k}=\dbigvee\limits_{j=1}^{n}a_{kj}\left( \underline{\delta _{i}}\right)
_{j}=\dbigvee\limits_{j=1}^{n}a_{kj}\delta _{ji}=a_{ki}=\left( T\underline{
\delta _{i}}\right) _{k}
\end{equation*}
for all $i,k=1,\ldots ,n$ and it follows that the action of $A_{T}$ is given
by $T$. The matrix $A_{T}$ is called the \emph{matrix corresponding to }$T$
in the canonical basis of $\mathcal{L}_{n}\left( \mathcal{B}\right) $. If $%
A= \left[ a_{ij}\right] _{n\times n}$ is a matrix on $\mathcal{L}_{n}\left( 
\mathcal{B}\right) $ then its transpose $\left[ a_{ji}\right] _{n\times n}$
is denoted by $A^{\ast }$.

\qquad It is straightforward to check that if $T:\mathcal{L}_{n}\left( 
\mathcal{B}\right) \longrightarrow \mathcal{L}_{m}\left( \mathcal{B}\right) $
and $S:\mathcal{L}_{m}\left( \mathcal{B}\right) \longrightarrow \mathcal{L}
_{k}\left( \mathcal{B}\right) $ then the matrix of $S\circ T$ is given by
the product $A_{S}A_{T}$, the matrix of $\lambda T$ for $\lambda \in 
\mathcal{B}$ is given by $\lambda A_{T}$ and if $S:\mathcal{L}_{n}\left( 
\mathcal{B}\right) \longrightarrow \mathcal{L}_{m}\left( \mathcal{B}\right) $
then the matrix of $S+T$ is $A_{S}+A_{T}$. Moreover, for all $\underline{a}
\in \mathcal{L}_{n}\left( \mathcal{B}\right) $ and $\underline{b}\in 
\mathcal{L}_{m}\left( \mathcal{B}\right) $ we check that $\left\langle T 
\underline{a},\underline{b}\right\rangle =\left\langle \underline{a},T^{\ast
}\underline{b}\right\rangle $ where $T^{\ast }:\mathcal{L}_{m}\left( 
\mathcal{B}\right) \longrightarrow \mathcal{L}_{n}\left( \mathcal{B}\right) $
is the linear map of matrix $A_{T}^{\ast }$ (and where we use the same
notation for the inner products on $\mathcal{L}_{n}\left( \mathcal{B}\right) 
$ and $\mathcal{L}_{m}\left( \mathcal{B}\right) $). Thus, linear maps always
have an adjoint. It is routine to check that the adjoint is unique. We thus
have, as with standard linear algebra, a natural isomorphism between the
*-algebra of linear maps and the *-algebra of Boolean matrices.

\qquad Invertibility of Boolean matrices was studied in \cite%
{Luce52,Rutherford63,Skornyakov86,Yoeli61} and the following result is
well-known. We present here a short proof which relies upon our previous
work with orthonormal bases and generalize the invertibility result to show
that invertible rectangular matrices have to be square. Note that if a
matrix $A$ is invertible, then its columns and its rows both form generating
families. We now show that these families are actually orthonormal bases of $%
\mathcal{L}_{n}\left( \mathcal{B}\right) $ and therefore are stochastic.

\begin{thm}
\label{Inverse}Let $A$ be an $n\times m$ Boolean matrix. The following are
equivalent:

\begin{enumerate}
\item $A$ is invertible, i.e. there exists a (necessarily unique) $m\times n$
Boolean matrix $A^{-1}$ such that $A^{-1}A=I_{n}$ and $AA^{-1}=I_{m}$,

\item $A$ is unitary,i.e. $n=m$ and $AA^{\ast }=A^{\ast }A=I_{n}$,

\item The columns of $A$ form an orthonormal basis for $\mathcal{L}
_{n}\left( \mathcal{B}\right) $,

\item The rows of $A$ form an orthonormal basis of $\mathcal{L}_{m}\left( 
\mathcal{B}\right) $.
\end{enumerate}

In particular, if any of 1-4 holds, then $n=m$.
\end{thm}

\begin{pf}
Assume (3) holds. Then by Theorem (\ref{Dimension}), there are $n$ columns
of $A$ and thus $n=m$. By Theorem (\ref{Duality}), the rows of $A$ are a
basis for $\mathcal{L}_{n}(\mathcal{B)}$ as well, so (4) holds. The same
reasoning shows that (4) implies (3) and in particular $n=m$ again.

Moreover, let us denote the columns of $A$ by $\underline{a_{1}},\ldots , 
\underline{a_{m}}$ and the rows of $A$ by $\underline{r_{1}},\ldots , 
\underline{r_{n}}$. By construction $A^{\ast }A=\left[ \left\langle 
\underline{a_{i}},\underline{a_{j}}\right\rangle \right] _{m\times m}$ and $%
AA^{\ast }=\left[ \left\langle \underline{r_{i}},\underline{r_{j}}
\right\rangle \right] _{n\times n}$ so $A$ is unitary if and only if both
(3) and (4) holds. Since (3) and (4) are equivalent and imply $n=m$, either
imply (2).

Assume now that $A$ is invertible and write $A=\left[ a_{ij}\right]
_{n\times m}$ and $A^{-1}=\left[ b_{ij}\right] _{m\times n}$. Then $%
A^{-1}A=I_{m}$ and $AA^{-1}=I_{n}$ implies that $\dbigvee
\limits_{j=1}^{m}a_{ij}=\dbigvee\limits_{j=1}^{n}b_{ij}=1$ and $%
b_{ki}a_{ij}=0$ ($k\not=j\in \left\{ 1,\ldots ,m\right\} $ and $i\in \left\{
1,\ldots n\right\} $) and $a_{ik}b_{kj}=0$ ($i\not=j\in \left\{ 1,\ldots
n\right\} $ and $k\in \left\{ 1,\ldots m\right\} $). Moreover if $i\in
\left\{ 1,\ldots ,n\right\} $ and $j\not=k\in \left\{ 1,\ldots ,m\right\} $
then: 
\begin{equation*}
a_{ij}a_{ik}=\left( \dbigvee\limits_{s=1}^{n}b_{ks}a_{ij}\right)
a_{ik}=\left( \dbigvee\limits_{\substack{ s=1  \\ s\not=i}}
^{n}a_{ij}b_{ks}\right) a_{ik}\leq \left(
\dbigvee\limits_{s=1,s\not=i}^{n}a_{ik}b_{ks}\right) =0\text{.}
\end{equation*}
Hence, the columns of $A^{\ast }$ form an orthonormal subset of $\mathcal{L}
_{m}\left( \mathcal{B}\right) $ and thus by Theorem (\ref{Duality}), the
columns of $A$ form an orthonormal basis of $\mathcal{L}_{n}\left( \mathcal{%
B }\right) $. So (1) implies (3) and the proof is complete.\qed
\end{pf}

\qquad As a consequence of Theorem (\ref{Inverse}), we see that invertible
operators are always isomorphisms by Lemma (\ref{IsometryChar}), since they
map the canonical basis to the orthonormal basis of their column vectors.

\bigskip Theorem (\ref{Inverse}) allows us to establish the following
remarkable fact: bases, as per Definition (\ref{Basis}), are necessarily
orthonormal, hence of cardinality the dimension of the Boolean vector space.
Thus, for Boolean vector spaces, being a basis in a traditional sense is the
same as being an orthonormal basis.

\begin{thm}
If $\mathcal{A=}\left\{ \underline{a_{1}},\ldots ,\underline{a_{m}}\right\} $
is a basis for $\mathcal{L}_{n}\left( \mathcal{B}\right) $ then $n=m$ and $%
\mathcal{A}$ is an orthonormal basis.
\end{thm}

\begin{pf}
Define 
\begin{equation*}
T:\left\vert 
\begin{array}{ccc}
\mathcal{L}_{n}\left( \mathcal{B}\right) & \longrightarrow & \mathcal{L}
_{m}\left( \mathcal{B}\right) \\ 
\underline{b} & \longmapsto & \left( b_{1},\ldots ,b_{m}\right)%
\end{array}
\right.
\end{equation*}
where $\underline{b}=\sum b_{i}\underline{a_{i}}$. Now $T$ is a linear
bijection. Denote the inverse of $T$ by $S$. It is easily checked that $S$
is a linear bijection and $ST=I_{n}$ and $TS=I_{m}$. We conclude from
Theorem (\ref{Inverse}) that $n=m$ and that matrix $A_{S}$ of $S$ is
unitary. It is easily checked that $S\underline{\delta _{i}}=\underline{%
a_{i} }$ for $i=1,\ldots ,n$, i.e. the columns of $A_{S}$ are the vectors $%
\underline{a_{1}},\ldots ,\underline{a_{n}}$ which by Theorem (\ref{Inverse}
) form an orthonormal basis.\qed
\end{pf}

\qquad We record the following observation as well:

\begin{cor}
Let $T:\mathcal{L}_{n}\left( \mathcal{B}\right) \longrightarrow \mathcal{L}
_{m}\left( \mathcal{B}\right) $ be a linear bijection. Then $n=m$ and $T$ is
an isomorphism.
\end{cor}

\bigskip \qquad In view of Theorem (\ref{Inverse}), we introduce a type of
matrix which will be of great interest to us in the next section. First,
given $A,B$ two $n\times n$ matrices, we shall say that $A\leq B$ when $%
\left\langle A\underline{a},\underline{b}\right\rangle \leq \left\langle B 
\underline{a},\underline{b}\right\rangle $ for all $\underline{a},\underline{
b}\in \mathcal{L}_{n}\left( \mathcal{B}\right) $. The relation $\leq $ is
easily seen to be an order on the set of $n\times n$ matrices. It is shown
in \cite{Luce52} that $\left[ a_{ij}\right] _{n\times n}\leq \left[ b_{ij} %
\right] _{n\times n}$ if and only $a_{ij}\leq b_{ij}$ for all $i,j\in
\left\{ 1,\ldots ,n\right\} $. Now we set:

\begin{defn}
A matrix $A$ is stochastic when $A^{\ast }A\geq I$ and $AA^{\ast }\leq I$.
\end{defn}

\qquad It is shown in \cite{Luce52} that products of stochastic matrices are
stochastic matrices, and that a matrix is stochastic if and only if it maps
stochastic vectors to stochastic vectors, or equivalently when its columns
are stochastic vectors.

\qquad Note that $A$ is unitary, or equivalently invertible, if and only if $%
A$ and $A^{\ast }$ are both stochastic. So unitarity is the same as
bi-stochasticity. As an interesting observation, if we call a matrix $A$
symmetric when $A^{\ast }=A$, then a symmetric stochastic matrix is always a
unitary of order 2, namely $A^{2}=I$. Conversely, if $A^{2}=I$ then $A$ is
invertible with $A^{-1}=A^{\ast }$, so symmetric stochastic matrices are
exactly given by unitaries of order 2, i.e. a reflection.

\qquad We have encountered such matrices before. Example (\ref{CyclicBasis})
shows how to obtain such reflections. Let $\underline{a}=\left( a_{1},\ldots
.a_{n}\right) $ be a stochastic vector. Then the matrix 
\begin{equation*}
A=\left[ 
\begin{array}{cccc}
a_{1} & a_{2} & \cdots & a_{n} \\ 
a_{2} & a_{3} & \cdots & a_{1} \\ 
\vdots & \vdots &  & \vdots \\ 
a_{n} & a_{1} & \cdots & a_{n-1}%
\end{array}
\right]
\end{equation*}
is symmetric and stochastic.

\qquad Note however that the product of reflections need not be a
reflection, as the product of the reflections $\left[ 
\begin{array}{ccc}
0 & 1 & 0 \\ 
1 & 0 & 0 \\ 
0 & 0 & 1%
\end{array}
\right] $ and $\left[ 
\begin{array}{ccc}
1 & 0 & 0 \\ 
0 & 0 & 1 \\ 
0 & 1 & 0%
\end{array}
\right] $ is given by $\left[ 
\begin{array}{ccc}
0 & 0 & 1 \\ 
1 & 0 & 0 \\ 
0 & 1 & 0%
\end{array}
\right] $ which is not a reflection.

\section{Invariant Vectors}

\qquad Eigenvalues and eigenvectors of Boolean matrices have been previously
studied \cite{Blyth67,Rutherford63b,Sindak75,Tan98}. Though invariant
vectors are special case of eigenvectors, as far as we know the results in
this section are new.

\qquad The following consequence of Lemma (\ref{Redux1}) will be used.

\begin{lem}
\label{Redux2}If $\underline{a}\in \mathcal{L}_{n}\left( \mathcal{B}\right) $
then there exists an orthovector $\underline{b}\in \mathcal{L}_{n}\left( 
\mathcal{B}\right) $ such that $\left\Vert \underline{b}\right\Vert
=\left\Vert \underline{a}\right\Vert $ and $\underline{b}\leq \underline{a}$.
\end{lem}

\begin{pf}
Apply Lemma (\ref{Redux1}) with the interval $\left[ 0,\left\Vert
a\right\Vert \right] $ of $\mathcal{B}$ in lieu of $\mathcal{B}$.\qed
\end{pf}

\bigskip Let $A_{1},\ldots ,A_{m}$ be matrices on $\mathcal{L}_{n}\left( 
\mathcal{B}\right) $ with $A_{k}=\left[ a_{ij}^{k}\right] _{n\times n}$ ($%
k=1,\ldots ,m$). The \emph{joint trace} of $A_{1},\ldots ,A_{m}$ is 
\begin{equation*}
\limfunc{tr}\left( A_{1},\ldots ,A_{m}\right)
=\dbigvee\limits_{i=1}^{n}a_{ii}^{1}a_{ii}^{2}\ldots a_{ii}^{m}\text{.}
\end{equation*}
In particular, the \emph{trace} of $\left[ a_{ij}\right] _{n\times n}$ is
given by $\limfunc{tr}\left( A\right) =\dbigvee\limits_{i=1}^{n}a_{ii}$. A
vector $\underline{b}$ is an \emph{invariant vector} for $A$ if $A\underline{
b}=\underline{b}$, and more generally a \emph{common invariant vector }of $%
A_{1},\ldots ,A_{m}$ if $A_{i}\underline{b}=\underline{b}$ for $i=1,\ldots
,m $.

\begin{lem}
Let $A,B$ be two matrices on $\mathcal{L}_{n}\left( \mathcal{B}\right) $.
Then

\begin{enumerate}
\item $\limfunc{tr}\left( AB\right) =\limfunc{tr}(BA)$,

\item If $B$ is invertible then $\limfunc{tr}\left( BAB^{\ast }\right) = 
\limfunc{tr}(A)$.
\end{enumerate}
\end{lem}

\begin{pf}
We compute 
\begin{equation*}
\limfunc{tr}(AB)=\dbigvee\limits_{i=1}^{n}\left( AB\right)
_{ii}=\dbigvee\limits_{i=1}^{n}\dbigvee\limits_{k=1}^{n}a_{ik}b_{ki}=
\dbigvee\limits_{k=1}^{n}\dbigvee\limits_{i=1}^{n}b_{ki}a_{ik}=\dbigvee
\limits_{k=1}^{n}\left( BA\right) _{kk}=\limfunc{tr}\left( BA\right) \text{.}
\end{equation*}

If $B$ is invertible then $B^{-1}=B^{\ast }$ by Theorem (\ref{Inverse}) and
thus (1) implies (2).\qed
\end{pf}

\begin{thm}
\label{StoInv}Stochastic matrices $A_{1},\ldots ,A_{m}$ on $\mathcal{L}
_{n}\left( \mathcal{B}\right) $ have a common invariant stochastic vector if
and only if $\limfunc{tr}\left( A_{1},\ldots ,A_{m}\right) =1$.
\end{thm}

\begin{pf}
Suppose $\underline{b}$ is a stochastic vector and $A_{i}\underline{b}=%
\underline{b}$ for $i=1,\ldots ,m$. Then $\dbigvee%
\limits_{j=1}^{n}a_{ij}^{k}b_{j}=b_{i}$ for $k=1,\ldots ,m$ and $i=1,\ldots
,n$. Multiplying both sides by $b_{i}$ and since $\underline{b}$ is
stochastic, we obtain $a_{ii}^{k}b_{i}=b_{i}$. Hence, $b_{i}\leq a_{ii}^{k}$%
, $k=1,\ldots ,m$, so $b_{i}\leq a_{ii}^{1}a_{ii}^{2}\ldots a_{ii}^{m}$.
Therefore 
\begin{equation*}
\limfunc{tr}\left( A_{1},\ldots ,A_{m}\right)
=\dbigvee\limits_{i=1}^{n}a_{ii}^{1}a_{ii}^{2}\ldots a_{ii}^{m}\geq
\dbigvee\limits_{i=1}^{n}b_{i}=1\text{.}
\end{equation*}%
Conversely, suppose $\limfunc{tr}\left( A_{1},\ldots ,A_{m}\right)
=\dbigvee\limits_{i=1}^{n}a_{ii}^{1}a_{ii}^{2}\ldots a_{ii}^{m}=1$. By
Lemma\ (\ref{Redux1}), there exists a stochastic vector $\underline{b}%
=\left( b_{1},\ldots ,b_{n}\right) $ such that $b_{j}\leq
a_{jj}^{1}a_{jj}^{2}\ldots a_{jj}^{m}$. Since $b_{j}\leq a_{jj}^{k}$ ($%
k=1,\ldots ,m$) and $A_{k}$ is stochastic, we have that $a_{ij}^{k}b_{j}=0$
for $i\not=j$, $i,j=1,\ldots ,n$ and $k=1,\ldots ,m$. Hence 
\begin{equation*}
\left( A_{k}\underline{b}\right)
_{i}=\dbigvee\limits_{j=1}^{n}a_{ij}^{k}b_{j}=a_{ii}^{k}b_{i}=b_{i}\text{.}
\end{equation*}%
Therefore, $A_{k}\underline{b}=\underline{b}$ ($k=1,\ldots ,m$) so $%
\underline{b}$ is a common invariant stochastic vector for $A_{1},\ldots
,A_{m}$.\qed
\end{pf}

\begin{cor}
\label{StoInv1}A stochastic matrix $A$ has an invariant stochastic vector if
and only if $\limfunc{tr}(A)=1$.
\end{cor}

\begin{cor}
If $A$ is a stochastic matrix and $B$ is invertible on $\mathcal{L}
_{n}\left( \mathcal{B}\right) $ then $A$ has an invariant stochastic vector
if and only if $BAB^{\ast }$ does.
\end{cor}

\begin{cor}
\label{StoInv3}A stochastic vector $\underline{b}=\left( b_{1},\ldots
,b_{n}\right) $ is a common invariant vector for stochastic matrices $%
A_{1},\ldots ,A_{m}$ if and only if $b_{i}\leq a_{ii}^{1}a_{ii}^{2}\ldots
a_{ii}^{m}$ for all $i=1,\ldots ,n$.
\end{cor}

\qquad Stochastic matrices $A_{1},\ldots ,A_{m}$ on $\mathcal{L}_{n}\left( 
\mathcal{B}\right) $ are \emph{simultaneously reducible} if there exists an
invertible matrix $B$ on $\mathcal{L}_{n}\left( \mathcal{B}\right) $ and
matrices $C_{1},\ldots ,C_{m}$ on $\mathcal{L}_{n-1}\left( \mathcal{B}
\right) $ such that for $i=1,\ldots ,m$ we have 
\begin{equation*}
A_{i}=B\left[ 
\begin{array}{cc}
1 & 0 \\ 
0 & C_{i}%
\end{array}
\right] B^{\ast }\text{.}
\end{equation*}
Notice that the matrices $C_{1},\ldots ,C_{m}$ are stochastic since $B^{\ast
}A_{i}B=\left[ 
\begin{array}{cc}
1 & 0 \\ 
0 & C_{i}%
\end{array}
\right] $. In particular, if there is only one matrix $A$ in the above
definition, we say that $A$ is \emph{reducible.}

\begin{thm}
\label{UnitReduce}Unitary matrices $A_{1},\ldots ,A_{m}$ on $\mathcal{\ L}
_{n}\left( \mathcal{B}\right) $ are simultaneously reducible if and only if $%
\limfunc{tr}\left( A_{1},\ldots ,A_{m}\right) =1$.
\end{thm}

\begin{pf}
If $A_{1},\ldots ,A_{m}$ are simultaneously reducible then $A_{i}=B\left[ 
\begin{array}{cc}
1 & 0 \\ 
0 & C_{i}%
\end{array}
\right] B^{\ast }$ for some invertible matrix $B$ and some matrix $C_{i}$, $%
i=1,\ldots ,m$. Since $B$ is unitary, $B\underline{\delta _{1}}$ is
stochastic and 
\begin{equation*}
A_{i}\left( B\underline{\delta _{1}}\right) =B\left[ 
\begin{array}{cc}
1 & 0 \\ 
0 & C_{i}%
\end{array}
\right] \underline{\delta _{1}}=B\underline{\delta _{1}}
\end{equation*}
for $i=1,\ldots ,m$. Hence, $A_{1},\ldots ,A_{m}$ have a common invariant
vector, and thus by Theorem (\ref{StoInv}) we have $\limfunc{tr}
(A_{1},\ldots ,A_{m})=1$.

Conversely, assume that $\limfunc{tr}\left( A_{1},\ldots ,A_{m}\right) =1$.
Then $A_{1},\ldots ,A_{m}$ have a common stochastic invariant vector $%
\underline{b}=\left( b_{1},\ldots ,b_{n}\right) $ by Theorem (\ref{StoInv}).
We define the symmetric stochastic matrix $B$ by 
\begin{equation*}
B=\left[ 
\begin{array}{ccccc}
b_{1} & b_{2} & b_{3} & \cdots & b_{n} \\ 
b_{2} & b_{2}^{c} & 0 & \cdots & 0 \\ 
b_{3} & 0 & b_{3}^{c} & \cdots & 0 \\ 
\vdots &  &  &  &  \\ 
b_{n} & 0 & 0 & \cdots & b_{n}^{c}%
\end{array}
\right] \text{.}
\end{equation*}
Let $D_{i}=BA_{i}B$ for $i=1,\ldots ,m$. With the notation $A_{k}=\left[
a_{ij}^{k}\right] _{n\times n}$, we compute the $(1,1)$ entry of $D_{i}$ as 
\begin{equation*}
\dbigvee\limits_{j=1}^{n}b_{1j}\left(
\dbigvee\limits_{r=1}^{n}a_{jr}^{i}b_{r1}\right)
=\dbigvee\limits_{j=1}^{n}b_{j}\left(
\dbigvee\limits_{r=1}^{n}a_{jr}^{i}b_{r}\right)
=\dbigvee\limits_{j=1}^{n}b_{j}b_{j}=1\text{.}
\end{equation*}
Since a product of unitary matrices is unitary, $D_{i}$ is a unitary matrix
and thus must have the form 
\begin{equation*}
D_{i}=\left[ 
\begin{array}{cc}
1 & 0 \\ 
0 & C_{i}%
\end{array}
\right]
\end{equation*}
for some matrix $C_{i}$ ($i=1,\ldots ,m$). Since $A_{i}=BD_{i}B$ for $%
i=1,\ldots ,m$, we are finished.\qed
\end{pf}

\begin{cor}
\label{UnitReduce1}A unitary matrix $A$ is reducible if and only if $%
\limfunc{tr}(A)=1$.
\end{cor}

\qquad We now give an example to show that Theorem (\ref{UnitReduce}) does
not hold for stochastic matrices. Consider the stochastic matrix $A=\left[ 
\begin{array}{cc}
1 & 1 \\ 
0 & 0%
\end{array}
\right] $. It is of trace $1$, yet if it were reducible then there exists a
unitary $B$ such that $A=B\left[ 
\begin{array}{cc}
1 & 0 \\ 
0 & 1%
\end{array}
\right] B^{\ast }=I$ which is a contradiction.

\qquad Notice if $A$ is unitary and $\underline{b}$ is an invariant vector
for $A$, then $\underline{b}$ is also an invariant vector for $A^{\ast }$.
Indeed, $A\underline{b}=\underline{b}$ implies that $A^{\ast }\underline{b}
=A^{\ast }A\underline{b}=\underline{b}$.

\bigskip \qquad We now give an example that motivates the next result. Let $%
A=\left[ a_{ij}\right] _{3\times 3}$ be a $3\times 3$ symmetric stochastic
matrix. We shall show that $A$ has an invariant stochastic vector and hence $%
A$ is reducible. Indeed, we have that 
\begin{eqnarray*}
a_{11}^{c}a_{22}^{c}a_{33}^{c} &=&\left( a_{12}\vee a_{13}\right) \left(
a_{12}\vee a_{32}\right) \left( a_{13}\vee a_{23}\right) \\
&=&\left( a_{12}\vee a_{13}\right) \left( a_{12}\vee a_{23}\right) \left(
a_{13}\vee a_{23}\right) \\
&=&\left( a_{12}a_{12}\vee a_{12}a_{23}\vee a_{13}a_{12}\vee
a_{13}a_{23}\right) \left( a_{13}\vee a_{23}\right) \\
&=&a_{12}\left( a_{13}\vee a_{23}\right) =0\text{.}
\end{eqnarray*}
Thus $\limfunc{tr}\left( A\right) =\left(
a_{11}^{c}a_{22}^{c}a_{33}^{c}\right) ^{c}=0^{c}=1$ so the result follows
from Corollaries (\ref{StoInv1}) and (\ref{UnitReduce1}). The next theorem
generalizes this calculation.

\begin{thm}
\label{OddInv}If $A$ is an $n\times n$ symmetric stochastic matrix with $n$
odd, then $A$ has an invariant stochastic vector.
\end{thm}

\begin{pf}
Since $A=\left[ a_{ij}\right] _{n\times n}$ is symmetric, we have that 
\begin{eqnarray*}
a_{11}^{c}a_{22}^{c}\ldots a_{nn}^{c} &=&\left( a_{12}\vee a_{13}\vee \ldots
\vee a_{1n}\right) \left( a_{12}\vee a_{23}\vee \ldots \vee a_{2n}\right) \\
&&\ldots \left( a_{1n}\vee a_{2n}\vee \ldots \vee a_{n-1,n}\right) \text{.}
\end{eqnarray*}
Since $A$ is stochastic, we conclude that if we expand the right hand-side,
the only nonzero terms are of the form $a_{ij}a_{ij}a_{rs}a_{rs}\ldots
a_{uv}a_{uv}$ with $i\not=r$, $r\not=u$ and so on. By construction, there
are $n$ factors in this product. This would imply that $n$ must be even.
This is a contradiction, so all terms in the expansion are zero and thus 
\begin{equation*}
\limfunc{tr}\left( A\right) =\left( a_{11}^{c}a_{22}^{c}\ldots
a_{nn}^{c}\right) ^{c}=1\text{.}
\end{equation*}

The result follows from Corollary (\ref{StoInv1}).\qed
\end{pf}

\bigskip \qquad We now show that Theorem\ (\ref{OddInv}) does not hold if $n$
is even. Consider the stochastic symmetric matrix $A=\left[ 
\begin{array}{cc}
0 & 1 \\ 
1 & 0%
\end{array}
\right] $. Then $\limfunc{tr}(A)=0$ so $A$ has no stochastic invariant
vector. Now, generalizing, we see that if $B$ is a $k\times k$ stochastic
symmetric matrix, then $\left[ 
\begin{array}{cc}
0 & B \\ 
B & 0%
\end{array}
\right] $ has trace $0$ and thus has no invariant stochastic vector. Thus,
for all even $n$ there exists a stochastic symmetric $n\times n$ matrix with
no invariant stochastic vector.

\bigskip \qquad We can find more invariant stochastic vectors in the natural
way. An \emph{invariant orthogonal set} for matrices $A_{1},\ldots ,A_{m}$
on $\mathcal{L}_{n}\left( \mathcal{B}\right) $ is a set of mutually
orthogonal invariant vectors for $A_{1},\ldots ,A_{m}$. For example, if $%
\underline{b},\underline{c}$ are stochastic vectors, then $\left\{ 
\underline{b},\underline{c}\right\} $ is an invariant orthogonal set for the
unitary matrix $A$ if and only if $c_{i}\leq a_{ii}b_{i}^{c}$ for $%
i=1,\ldots ,n$ or equivalently $b_{i}\leq a_{ii}c_{i}^{c}$ for $i=1,\ldots
,n $.

\begin{thm}
\label{UnitaryReduxMultiple}A unitary matrix $A$ possesses an invariant
orthogonal set of $m$ stochastic vectors if and only if there exists an
invertible matrix $B$ such that 
\begin{equation*}
A=B\left[ 
\begin{array}{cc}
I_{m} & 0 \\ 
0 & C%
\end{array}
\right] B^{\ast }
\end{equation*}
where $I_{m}$ is the identity operator on $\mathcal{L}_{m}\left( \mathcal{B}
\right) $.
\end{thm}

\begin{pf}
Suppose $A$ is an $n\times n$ matrix with the given form. Then $m\leq n$ and
we can define $\underline{b_{j}}=B\underline{\delta _{j}}$, $j=1,\ldots ,m$.
We conclude from Theorem (\ref{Inverse}) that $\underline{b_{1}},\ldots , 
\underline{b_{m}}$ are stochastic vectors and we have $A\underline{b_{j}}= 
\underline{b_{j}}$ for $j=1,\ldots ,m$ by construction. Moreover, for $%
i\not=j$ we have 
\begin{equation*}
\left\langle \underline{b_{j}},\underline{b_{i}}\right\rangle =\left\langle
B \underline{\delta _{i}},B\underline{\delta _{j}}\right\rangle
=\left\langle B^{\ast }B\underline{\delta _{i}},\underline{\delta _{j}}%
\right\rangle =\left\langle \underline{\delta _{i}},\underline{\delta _{j}}%
\right\rangle =0 \text{.}
\end{equation*}

Hence $\left\{ \underline{b_{1}},\ldots ,\underline{b_{m}}\right\} $ is an
invariant orthogonal set of stochastic vectors.

\qquad Conversely, suppose that $A$ possesses an invariant orthogonal set of
stochastic vectors $\left\{ \underline{b_{1}},\ldots ,\underline{b_{m}}%
\right\} $ and write $\underline{b_{j}}=\left( b_{1j},\ldots ,b_{nj}\right) $
for $j=1,\ldots ,m$. Letting 
\begin{equation*}
B_{1}=\left[ 
\begin{array}{ccccc}
b_{11} & b_{21} & b_{31} & \cdots & b_{n1} \\ 
b_{21} & b_{21}^{c} & 0 & \cdots & 0 \\ 
b_{31} & 0 & b_{31}^{c} & \cdots & 0 \\ 
\vdots &  &  &  & \vdots \\ 
b_{n1} & 0 & \cdots & 0 & b_{n1}^{c}%
\end{array}%
\right]
\end{equation*}%
and $D_{1}=B_{1}AB_{1}$ as in the proof of Theorem (\ref{UnitReduce}), we
have that 
\begin{equation*}
D_{1}=\left[ 
\begin{array}{cc}
1 & 0 \\ 
0 & C_{1}%
\end{array}%
\right]
\end{equation*}%
where $C_{1}$ is a stochastic matrix and $A=B_{1}D_{1}B_{1}$. Letting $C_{1}=%
\left[ c_{ij}\right] _{(n-1)\times (n-1)}$ and $D_{1}=\left[ d_{ij}\right]
_{n\times n}$ we have 
\begin{eqnarray*}
c_{11} &=&d_{22}=\dbigvee\limits_{j=1}^{2}b_{2j}\left(
\dbigvee\limits_{k=1}^{2}a_{jk}b_{k2}\right) \\
&=&b_{21}\left( a_{11}b_{21}\vee a_{12}b_{21}^{c}\right) \vee
b_{21}^{c}\left( a_{21}b_{21}\vee a_{22}b_{21}^{c}\right) \\
&=&a_{11}b_{21}\vee a_{22}b_{21}^{c}\text{.}
\end{eqnarray*}%
More generally 
\begin{equation*}
c_{ii}=d_{i+1,i+1}=a_{ii}b_{i+1,1}\vee a_{i+1,i+1}b_{i+1,1}^{c}
\end{equation*}%
for $i=1,\ldots ,n-1$. Hence 
\begin{equation*}
\limfunc{tr}\left( C_{1}\right) =\dbigvee\limits_{i=1}^{n-1}\left(
a_{ii}b_{i+1,1}\vee a_{i+1,i+1}b_{i+1,1}^{c}\right)
=\dbigvee_{i=1}^{n}a_{ii}b_{i,1}^{c}\text{.}
\end{equation*}%
Since $b_{i2}\leq a_{ii}b_{i1}^{c}$ ($i=1,\ldots ,n$), we conclude that $%
\underline{b_{2}}$ is an invariant stochastic vector of $C_{1}$ by Corollary
(\ref{StoInv3}). Hence, there exists a symmetric stochastic matrix $B_{2}$
such that 
\begin{equation*}
C_{1}=B_{2}\left[ 
\begin{array}{cc}
1 & 0 \\ 
0 & C_{2}%
\end{array}%
\right] B_{2}\text{.}
\end{equation*}%
It follows that 
\begin{eqnarray*}
A &=&B_{1}\left[ 
\begin{array}{cc}
1 & 0 \\ 
0 & B_{2}\left[ 
\begin{array}{cc}
1 & 0 \\ 
0 & C_{2}%
\end{array}%
\right] B_{2}%
\end{array}%
\right] B_{1} \\
&=&B_{1}\left[ 
\begin{array}{cc}
1 & 0 \\ 
0 & B_{2}%
\end{array}%
\right] \left[ 
\begin{array}{ccc}
1 & 0 & 0 \\ 
0 & 1 & 0 \\ 
0 & 0 & C_{2}%
\end{array}%
\right] \left[ 
\begin{array}{cc}
1 & 0 \\ 
0 & B_{2}%
\end{array}%
\right] B_{1} \\
&=&B_{3}\left[ 
\begin{array}{cc}
I_{2} & 0 \\ 
0 & C_{2}%
\end{array}%
\right] B_{3}^{\ast }
\end{eqnarray*}%
with $B_{3}=B_{1}\left[ 
\begin{array}{cc}
1 & 0 \\ 
0 & B_{2}%
\end{array}%
\right] $. The proof is then completed by a simple induction.\qed
\end{pf}

\qquad Theorem (\ref{UnitaryReduxMultiple}) can be easily generalized to the
following:

\begin{cor}
Unitary matrices $A_{1},\ldots ,A_{m}$ possess an invariant orthogonal set
of stochastic vectors if and only if there exists an invertible matrix $B$
and matrices $C_{1},\ldots ,C_{n}$ such that 
\begin{equation*}
A_{i}=B\left[ 
\begin{array}{cc}
I_{m} & 0 \\ 
0 & C_{i}%
\end{array}
\right] B^{\ast }
\end{equation*}
for $i=1,\ldots ,m$ and $I_{m}$ the identity operator on $\mathcal{L}
_{m}\left( \mathcal{B}\right) $.
\end{cor}

\qquad We now illustrate Theorem (\ref{UnitaryReduxMultiple}) with an
example. Let $\mathcal{B}$ be the power set of $\left\{ 1,2,3,4,5\right\}
=\Omega $ endowed with its natural Boolean algebra structure. Consider the
stochastic symmetric matrix $A$ over $\mathcal{L}_{5}\left( \mathcal{B}
\right) $ defined by 
\begin{equation*}
A=\left[ 
\begin{array}{ccccc}
\left\{ 1\right\} & \left\{ 2\right\} & \left\{ 3\right\} & \left\{ 4\right\}
& \left\{ 5\right\} \\ 
\left\{ 2\right\} & \left\{ 4,5\right\} & \emptyset & \emptyset & \left\{
1,3\right\} \\ 
\left\{ 3\right\} & \emptyset & \left\{ 4,5\right\} & \left\{ 1\right\} & 
\{2\} \\ 
\left\{ 4\right\} & \emptyset & \left\{ 1\right\} & \left\{ 2,3,5\right\} & 
\emptyset \\ 
\left\{ 5\right\} & \left\{ 1,3\right\} & \left\{ 2\right\} & \emptyset & 
\left\{ 4\right\}%
\end{array}
\right] \text{.}
\end{equation*}
There are many stochastic invariant vectors for $A$ and we choose 
\begin{equation*}
\underline{b}=\left( \left\{ 1\right\} ,\emptyset ,\emptyset ,\left\{
2,3,5\right\} ,\left\{ 4\right\} \right) \text{.}
\end{equation*}
We now form the stochastic symmetric matrix 
\begin{equation*}
B=\left[ 
\begin{array}{ccccc}
\left\{ 1\right\} & \emptyset & \emptyset & \left\{ 2,3,5\right\} & \left\{
4\right\} \\ 
\emptyset & \Omega & \emptyset & \emptyset & \emptyset \\ 
\emptyset & \emptyset & \Omega & \emptyset & \emptyset \\ 
\left\{ 2,3,5\right\} & \emptyset & \emptyset & \left\{ 1,4\right\} & 
\emptyset \\ 
\left\{ 4\right\} & \emptyset & \emptyset & \emptyset & \left\{
1,2,3,5\right\}%
\end{array}
\right]
\end{equation*}
We can then reduce $A$ by 
\begin{equation*}
BAB=\left[ 
\begin{array}{ccccc}
\Omega & \emptyset & \emptyset & \emptyset & \emptyset \\ 
\emptyset & \left\{ 4,5\right\} & \emptyset & \left\{ 2\right\} & \left\{
1,3\right\} \\ 
\emptyset & \emptyset & \left\{ 4,5\right\} & \left\{ 1,3\right\} & \left\{
2\right\} \\ 
\emptyset & \left\{ 2\right\} & \left\{ 1,3\right\} & \emptyset & \left\{
4,5\right\} \\ 
\emptyset & \left\{ 1,3\right\} & \left\{ 2\right\} & \left\{ 4,5\right\} & 
\emptyset%
\end{array}
\right] \text{.}
\end{equation*}
Thus 
\begin{equation*}
A=B\left[ 
\begin{array}{cc}
1 & 0 \\ 
0 & C%
\end{array}
\right] B
\end{equation*}
yet $\limfunc{tr}(C)=\left\{ 4,5\right\} \not=\Omega $ so no further
reduction is possible.

\section{Powers of Stochastic Matrices}

\qquad As mentioned in section 2, powers of stochastic matrices may be
important for the study of Boolean Markov chains. Various applications of
powers of lattice matrices are discussed in \cite{Cechlarova03,Tan01}. If $A$
is a Boolean matrix, the smallest natural number $p$ such that there exists
a natural number $e$ with $A^{e+p}=A^{e}$ is called the \emph{period} of $A$
and is denoted by $p(A)$. The smallest natural number $e$ such that $%
A^{e+p(A)}=A^{e}$ is called the \emph{exponent} or \emph{index} of $A$ and
is denoted by $e(A)$. It is known that for any $n\times n$ Boolean matrix $A$
, both $p(A)$ and $e(A)$ exist and $e(A)\leq \left( n-1\right) ^{2}+1$ \cite%
{Cechlarova03,Tan01}. We shall use:

\begin{defn}
Let $n\in \mathbb{N}$. The least common multiple of $\left\{ 1,2,\ldots
,n\right\} $ is denoted by $[n]$.
\end{defn}

\qquad It is also known that $p(A)$ divides $[n]$.

\qquad In this section, we show that for a stochastic matrix, we can improve
the upper bound for $e(A)$ to $e(A)\leq n-1$. Although we do not improve on $%
p(A)|[n]$, we give an alternative proof of this result for stochastic
matrices because it is embedded in our proof that $e(A)\leq n-1$.

\qquad If $A$ is a $2\times 2$ matrix, then it follows from the previous
known results that $A^{4}=A^{2}$. Moreover, it is easy to check that if $A$
is a $2\times 2$ stochastic matrix then $A^{3}=A$. In the same way, for $%
3\times 3$ matrix $A$ we have $A^{11}=A^{5}$. However, one can check that if 
$A$ is a $3\times 3$ stochastic matrix then $A^{8}=A^{2}$. Displaying the
first eight powers of $A$ would be cumbersome, so we refrain from doing so.
However, we can easily prove the special case that $A^{6}=I$ for any unitary 
$3\times 3$ matrix $A$. In this case, we have 
\begin{equation*}
A=\left[ 
\begin{array}{ccc}
a_{1} & b_{1} & c_{1} \\ 
a_{2} & b_{2} & c_{2} \\ 
a_{3} & b_{3} & c_{3}%
\end{array}
\right]
\end{equation*}
where each row and column is a stochastic vector. We then have 
\begin{eqnarray*}
A^{2} &=&\left[ 
\begin{array}{ccc}
a_{1}\vee a_{2}b_{1}\vee a_{3}c_{1} & b_{3}c_{1} & b_{1}c_{2} \\ 
a_{3}c_{2} & a_{2}b_{1}\vee b_{2}\vee b_{3}c_{2} & a_{2}c_{1} \\ 
b_{3}a_{2} & a_{3}b_{1} & a_{3}c_{1}\vee b_{3}c_{2}\vee c_{3}%
\end{array}
\right] , \\
A^{3} &=&\left[ 
\begin{array}{ccc}
a_{1}\vee a_{3}b_{1}c_{2}\vee a_{2}b_{3}c_{1} & a_{2}b_{1} & a_{3}c_{1} \\ 
a_{2}b_{1} & a_{2}b_{3}c_{1}\vee b_{2}\vee a_{3}b_{1}c_{2} & b_{3}c_{2} \\ 
a_{3}c_{1} & b_{3}c_{2} & a_{3}b_{1}c_{2}\vee a_{2}b_{3}c_{1}\vee c_{3}%
\end{array}
\right] \text{.}
\end{eqnarray*}
Since $A^{3}$ is symmetric and unitary (as a product of unitary, or by
inspection), we conclude that $A^{6}=A^{3}A^{3}=I$.

\qquad From these observations and our work in Section 5, we can already
draw some interesting conclusions. For example, let $A$ be a $3\times 3$
unitary matrix with $\limfunc{tr}(A)=1$. Applying Corollary (\ref%
{UnitReduce1}), there exists an invertible matrix $B$ and a $2\times 2$
unitary matrix $C$ such that 
\begin{equation}
A=B\left[ 
\begin{array}{cc}
1 & 0 \\ 
0 & C%
\end{array}
\right] B^{\ast }\text{.}  \label{Eq61}
\end{equation}

Since $C$ is symmetric (all $2\times 2$ unitaries are), we have $C^{2}=I$
and thus 
\begin{equation*}
A^{2}=B\left[ 
\begin{array}{cc}
1 & 0 \\ 
0 & C^{2}%
\end{array}
\right] B^{\ast }=I\text{.}
\end{equation*}
We conclude that any $3\times 3$ unitary matrix $A$ with $\limfunc{tr}(A)=1$
is symmetric.

\qquad As another example, let $A$ be a $4\times 4$ unitary matrix with $%
\limfunc{tr}(A)=1$. As before, there exists an invertible matrix $B$ such
that (\ref{Eq61}) holds where $C$ is now a $3\times 3$ unitary matrix. Since 
$C^{6}=I$, we conclude that $A^{6}=I$ and thus $A^{3}$ is symmetric.

\bigskip \qquad We now begin the proof of the main result of this section.
Let $A=\left[ a_{ij}\right] _{n\times n}$ be a stochastic matrix on $%
\mathcal{L}_{n}\left( \mathcal{B}\right) $. We shall use:

\begin{defn}
A nonzero element of $\mathcal{B}$ of the form 
\begin{equation*}
a_{i_{1}1}a_{i_{2}2}\ldots a_{i_{n}n}
\end{equation*}
for $i_{1},\ldots ,i_{n}\in \left\{ 1,\ldots ,n\right\} $ is called an atom
of $A$.
\end{defn}

\qquad Of course there are a finite numbers of atoms of $A$.

\begin{lem}
\label{Atoms}Let $A=\left[ a_{ij}\right] _{n\times n}$ be a stochastic
matrix on $\mathcal{L}_{n}\left( \mathcal{B}\right) $. Let $\omega
_{1},\ldots ,\omega _{m}$ be the distinct atoms of $A$.

\begin{enumerate}
\item If $i,j\in \left\{ 1,\ldots ,m\right\} $ and $i\not=j$ then $\omega
_{i}\omega _{j}=0$,

\item $\dbigvee\limits_{i=1}^{m}\omega _{i}=1$,

\item For all $i,j\in \left\{ 1,\ldots ,n\right\} $ we have $a_{ij}=\dbigvee
\left\{ \omega _{k}:\omega _{k}\leq a_{ij}\right\} $,

\item If $\omega _{i}\leq a_{kj}$ then $A\omega _{i}\underline{\delta _{j}}
=\omega _{i}\underline{\delta _{k}}$.
\end{enumerate}
\end{lem}

\begin{pf}
For (1),\ letting $\omega _{i}=a_{i_{1}1}a_{i_{2}2}\ldots a_{i_{n}n}$ and $%
\omega _{j}=a_{j_{1}1}a_{j_{2}2}\ldots a_{j_{n}n}$, if $i\not=j$ then $%
i_{k}\not=j_{k}$ for some $k\in \left\{ 1.\ldots ,n\right\} $ and thus $%
\omega _{j}\omega _{i}=0$ since $a_{i_{k}k}a_{j_{k}k}=0$.

(2) will follow from (3). For (3), since 
\begin{equation*}
a_{11}=\dbigvee \left\{ a_{11}\left( a_{i_{2}2}\ldots a_{i_{n}n}\right)
:i_{2},\ldots ,i_{n}=1,\ldots ,n\right\}
\end{equation*}
as $A$ is stochastic, the results holds for $a_{11}$. It holds similarly for 
$a_{ij}$ with $i,j\in \left\{ 1,\ldots ,n\right\} $. Last, for (4), if $%
\omega _{i}\leq a_{kj}$ then 
\begin{eqnarray*}
A\omega _{i}\underline{\delta _{j}} &=&\omega _{i}A\underline{\delta _{j}}
=\omega _{i}\left( a_{1j},a_{2j},\ldots ,a_{nj}\right) =\left( \omega
_{i}a_{1j},\ldots ,\omega _{i}a_{nj}\right) \\
&=&\omega _{i}a_{jk}\underline{\delta _{k}}=\omega _{i}\underline{\delta
_{k} }\text{.}
\end{eqnarray*}
This concludes our proof.\qed
\end{pf}

\bigskip \qquad The main result for this section is:

\begin{thm}
If $A$ is a stochastic $n\times n$ matrix then $A^{[n]+n-1}=A^{n-1}$.
\end{thm}

\begin{pf}
Let $\omega _{1},\ldots ,\omega _{m}$ be the distinct atoms of $A$. By Lemma
(\ref{Atoms},2), we have $\underline{\delta _{i}}=\sum_{j=1}^{m}\omega _{j} 
\underline{\delta _{i}}$ for all $i\in \left\{ 1,\ldots n\right\} $. Since $%
\left\{ \underline{\delta _{1}},\ldots ,\underline{\delta _{n}}\right\} $ is
a basis for $\mathcal{L}_{n}\left( \mathcal{B}\right) $, the set $\left\{
\omega _{j}\underline{\delta _{i}}:i=1,\ldots n;j=1,\ldots m\right\} $ is a
generating set of $\mathcal{L}_{n}\left( \mathcal{B}\right) $. Set $%
r=[n]+n-1 $. If we can show that $A^{r}\omega _{j}\underline{\delta _{i}}
=A^{n-1}\omega _{j}\underline{\delta _{i}}$ for $i=1,\ldots n$ and $%
j=1,\ldots m$ then we are done.

Consider first $\omega _{1}\underline{\delta _{1}}$ and call the vectors $%
A^{0}\omega _{1}\underline{\delta _{1}}$, $A\omega _{1}\underline{\delta
_{1} }$, $A^{2}\omega _{1}\underline{\delta _{1}},\ldots ,A^{n-1}\omega _{1} 
\underline{\delta _{1}}$ the \emph{iterates} of $A$ \emph{at} $\omega _{1} 
\underline{\delta _{1}}$. By Lemma (\ref{Atoms},4), the iterates of $\omega
_{1}\underline{\delta _{1}}$ have the form: $\omega _{1}\underline{\delta
_{1}},$ $\omega _{1}\underline{\delta _{i_{1}}},\omega _{1}\underline{\delta
_{i_{2}}},\ldots ,\omega _{1}\underline{\delta _{i_{n-1}}}$ for $%
i_{1},\ldots ,i_{n-1}\in \left\{ 1,\ldots ,n\right\} $.

Suppose there is only one distinct iterate of $A$ at $\omega _{1}\underline{
\delta _{1}}$. Then 
\begin{equation*}
A\omega _{1}\underline{\delta _{1}}=\omega _{1}\underline{\delta _{i_{1}}}
=\omega _{1}\underline{\delta _{1}}\text{.}
\end{equation*}
Then we have 
\begin{equation}
A^{n-1}\omega _{1}\underline{\delta _{1}}=A^{n}\omega _{1}\underline{\delta
_{1}}=\ldots =A^{r}\omega _{1}\underline{\delta _{1}}\text{.}  \label{Eq62}
\end{equation}

Suppose now there are two distinct iterates of $A$ at $\omega _{1}\underline{%
\delta _{1}}$. Then $\omega _{1}\underline{\delta _{1}}\not=\omega _{1}%
\underline{\delta _{i_{1}}}$. If $\omega _{1}\underline{\delta _{i_{2}}}%
=\omega _{1}\underline{\delta _{i_{1}}}$ then 
\begin{equation*}
\omega _{1}\underline{\delta _{i_{3}}}=A\omega _{1}\underline{\delta _{i_{2}}%
}=A\omega _{1}\underline{\delta _{i_{1}}}=\omega _{1}\underline{\delta
_{i_{2}}}=\omega _{1}\underline{\delta _{i_{1}}}
\end{equation*}%
and we can conclude again that (\ref{Eq62}) holds. Otherwise, $A^{2}\omega
_{1}\underline{\delta _{1}}=\omega _{1}\underline{\delta _{1}}$ and thus $%
A^{n-1}\omega _{1}\underline{\delta _{1}}=\omega _{1}\underline{\delta _{1}}$
or $A^{n-1}\omega _{1}\underline{\delta _{1}}=\omega _{1}\underline{\delta
_{i_{1}}}$. Either way, we have 
\begin{equation}
A^{2+(n-1)}\omega _{1}\underline{\delta _{1}}=A^{n-1}\omega _{1}\underline{%
\delta _{1}}\text{.}  \label{Eq63}
\end{equation}

Suppose instead that there are three distinct iterates of $A$ at $\omega
_{1} \underline{\delta _{1}}$. Thus $\omega _{1}\underline{\delta _{1}},$ $%
\omega _{1}\underline{\delta _{i_{1}}}$ and $\omega _{1}\underline{\delta
_{i_{2}}}$ are distinct. If $\omega _{1}\underline{\delta _{i_{3}}}=\omega
_{1} \underline{\delta _{i_{2}}}$ then $A^{r}\omega _{1}\underline{\delta
_{1}} =A^{n-1}\omega _{1}\underline{\delta _{1}}=\omega _{1}\underline{%
\delta _{i_{3}}}$ so (\ref{Eq62}) holds again. If $\omega _{1}\underline{%
\delta _{i_{1}}}=\omega _{1}\underline{\delta _{i_{3}}}$ then $A\omega _{1} 
\underline{\delta _{1}}\in \left\{ \omega _{1}\underline{\delta _{1}},\omega
_{1}\underline{\delta _{i_{2}}}\right\} $ and (\ref{Eq63}) holds. If $\omega
_{1}\underline{\delta _{1}}=\omega _{1}\underline{\delta _{i_{3}}}$ then $%
A^{n-1}\omega _{1}\underline{\delta _{1}}\in \left\{ \omega _{1}\underline{
\delta _{1}},\omega _{1}\underline{\delta _{i_{1}}},\omega _{1}\underline{
\delta _{i_{2}}}\right\} $ and we have 
\begin{equation}
A^{3+n-1}\omega _{1}\underline{\delta _{1}}=A^{n-1}\omega _{1}\underline{
\delta _{1}}\text{.}  \label{Eq64}
\end{equation}

Generalizing this observation, suppose that all the iterates $\omega _{1} 
\underline{\delta _{1}},\omega _{1}\underline{\delta _{i_{1}}},\ldots ,$ $%
\omega _{1}\underline{\delta _{i_{n-1}}}$ are distinct. Since there are only 
$n$ possibilities for $A^{n}\omega _{1}\underline{\delta _{1}}$, we conclude
that $A^{n}\omega _{1}\underline{\delta _{1}}=\omega _{1}\underline{\delta
_{1}}$ or $\omega _{1}\underline{\delta _{i_{j}}}$ for some $j\in \left\{
1,\ldots ,n-1\right\} $. But then 
\begin{equation}
A^{t+\left( n-1\right) }\omega _{1}\underline{\delta _{1}}=A^{n-1}\omega
_{1} \underline{\delta _{1}}  \label{Eq65}
\end{equation}

for some $t\in \left\{ 1,2,\ldots ,n\right\} $. Notice (\ref{Eq63}) and (\ref%
{Eq64}) are special cases of (\ref{Eq65}).

Let us now suppose (\ref{Eq65}) holds for some $t\in \left\{ 1,\ldots
,n\right\} $. Since $r=kt+(n-1)$ for some $k\in \mathbb{N}$ we have 
\begin{eqnarray*}
A^{r}\omega _{1}\underline{\delta _{1}} &=&A^{kt+n-1}\omega _{1}\underline{
\delta _{1}}=\left( A^{t}\right) ^{k}A^{n-1}\omega _{1}\underline{\delta
_{1} } \\
&=&\left( A^{t}\right) ^{k-1}A^{t}A^{n-1}\omega _{1}\underline{\delta _{1}}
=\left( A^{t}\right) ^{k-1}A^{n-1}\omega _{1}\underline{\delta _{1}} \\
&=&\left( A^{t}\right) ^{k-2}A^{t}A^{n-1}\omega _{1}\underline{\delta _{1}}
=\left( A^{t}\right) ^{k-2}A^{n-1}\omega _{1}\underline{\delta _{1}} \\
&=&\ldots =A^{n-1}\omega _{1}\underline{\delta _{1}}\text{.}
\end{eqnarray*}
In a similar way, we can prove that $A^{r}\omega _{j}\underline{\delta _{i}}
=A^{n-1}\omega _{j}\underline{\delta _{i}}$ for $j=1,\ldots ,m$ and $%
i=1,\ldots ,n$, so the proof is complete.\qed
\end{pf}

\begin{cor}
If $A$ is an $n\times n$ unitary matrix then $A^{\left[ n\right] }=I$.
\end{cor}

\qquad As examples, $A^{15}=A^{3}$ for any $4\times 4$ stochastic matrix and 
$A^{64}=A^{4}$ for any $5\times 5$ stochastic matrix. We now give a final
example. Let $\left( a,b,c\right) $ be a stochastic vector and form the
stochastic matrix 
\begin{equation*}
A=\left[ 
\begin{array}{ccc}
b\vee c & a & 0 \\ 
a & b & a \\ 
0 & c & b\vee c%
\end{array}
\right] \text{.}
\end{equation*}
We then have 
\begin{equation*}
A^{2}=\left[ 
\begin{array}{ccc}
1 & 0 & a \\ 
0 & a\vee b & 0 \\ 
0 & c & c\vee b%
\end{array}
\right]
\end{equation*}
and $A^{2n+1}=A$, $A^{2n}=A^{2}$ for $n\in \mathbb{N}$. This example
illustrates an important difference between Boolean Markov chains and
traditional Markov chains given by real stochastic matrices. An important
property of traditional Markov chains is that the sites (called states in
the traditional case) can be decomposed into equivalence classes. This is
important because sites in the same equivalence class share a similar
behavior \cite{Giveon64}.

\qquad To be precise, let $M=\left[ p_{ij}\right] _{n\times n}$ be a real
stochastic matrix, i.e. $p_{ij}\geq 0$ and $\sum_{i=1}^{n}p_{ij}=1$ for
every $j=1,\ldots ,n$. The real $p_{ij}$ represents the transition
probability from site $j$ to site $i$. A site $i$ is \emph{accessible} from
a site $j$ if there exists $n\in \mathbb{N}$ such that $\left( M^{n}\right)
_{ij}>0$, and we then denote $j\rightarrow i$. It is easy to check that $%
\rightarrow $ is transitive and that the relation $\longleftrightarrow $
defined by $i\longleftrightarrow j\iff \left( i\rightarrow j\wedge
j\rightarrow i\right) $ is an equivalence relation on the sites of the
Markov chain.

\qquad Let us now extend this concept to Boolean Markov chains whose
transition matrix is a Boolean stochastic matrix $A$. Thus, $j\rightarrow i$
whenever $\left( A^{n}\right) _{ij}>0$ for some $n\in \mathbb{N}$. For the
example above, we note that $1\rightarrow 2$ and $2\rightarrow 3$ yet $%
1\not\rightarrow 3$. Thus $\rightarrow $ is not transitive. If we define $%
\longleftrightarrow $ by $i\longleftrightarrow j\iff \left( i\rightarrow
j\wedge j\rightarrow i\right) $ then we have, in the above example, that in
fact $1\longleftrightarrow 2$ and $2\longleftrightarrow 3$ yet $%
1\nleftrightarrow 3$. Hence $\longleftrightarrow $ is no longer an
equivalence relation.

\bibliographystyle{amsplain}
\bibliography{thesis}

\end{document}